\documentclass[a4paper, reqno, 11pt]{amsart}

\usepackage{amsfonts, amsmath, amssymb, amsthm, amscd, setspace, subfigure,mathtools,stmaryrd}
\usepackage{tikz-cd}
\usepackage[heightrounded,textwidth=400pt,top=1in,left=1in,right=1in,bottom=1in]{geometry}
\usepackage{layout}
\usepackage[utf8]{inputenc}
\usepackage{cite}
\usepackage{tikz}
\usepackage{hyperref}
\usetikzlibrary{calc}

\numberwithin{equation}{section}

\newcommand\be{\begin{eqnarray}}
\newcommand\ee{\end{eqnarray}}

\newcommand{\R}{{\mathbb R}}
\newcommand{\Z}{{\mathbb Z}}

\theoremstyle{definition}

\newtheorem{theorem}{Theorem}

\author[Krasnov]{Kirill Krasnov}
\email{kirill.krasnov@nottingham.ac.uk, ORCID: 0000-0003-2800-3767}
\address{School of Mathematical Sciences, University of Nottingham, Nottingham, NG7 2RD, UK}

\author[Shaw]{Adam Shaw}
\email{adam.shaw@nottingham.ac.uk}
\address{School of Mathematical Sciences, University of Nottingham, Nottingham, NG7 2RD, UK}

\title{Pleba\'nski complex}

\date{February 2025}

\begin{document}

\begin{abstract}\noindent As is very well-known, linearisation of the instanton equations on a 4-manifold gives rise to an elliptic complex of differential operators, the truncated (twisted) Hodge complex $\Lambda^0(\mathfrak{g}) \to \Lambda^1(\mathfrak{g})\to \Lambda^2_+(\mathfrak{g})$. Moreover, the linearisation of the full YM equations also fits into this framework, as it is given by the second map followed by its adjoint. We define and study properties of what we call the Pleba\'nski complex. This is a differential complex that arises by linearisation of the equations implying that a Riemannian 4-manifold is hyper-K\"ahler. We recall that these are most naturally stated as the condition that there exists a perfect $\Sigma^i\wedge \Sigma^j\sim\delta^{ij}$ triple $\Sigma^i, i=1,2,3$ of 2-forms that are closed $d\Sigma^i=0$. The Riemannian metric is encoded by the 2-forms $\Sigma^i$. We show that what results is an elliptic differential complex $TM \to S\to E\times \Lambda^1 \to E$, where $S$ is the tangent space to the space of perfect triples, and $E=\R^3$. We also show that, as in the case with instanton equations, the full Einstein equations $Ric=0$ also fit into this framework, their linearisation being given by the second map followed by its adjoint. Our second result concerns the elliptic operator that the Pleba\'nski complex defines. In the case of the instanton complex, operators appearing in the complex supplemented with their adjoints assemble to give the Dirac operator. We show how the same holds true for the Pleba\'nski complex. Supplemented by suitable adjoints, operators assemble into an elliptic operator that squares to the Laplacian and is given by the direct sum of two Dirac operators. 
\end{abstract}


\maketitle
\tableofcontents

\section{Introduction}

Since the pioneering works by Atiyah, Hitchin, Singer \cite{Atiyah:1978wi} and Donaldson \cite{Donaldson:1983wm} Yang-Mills (YM) instanton equations $F(A)_+=0$ on a 4-manifold played fundamental role in differential geometry. Linearisation of these equations around an instanton gives rise to an elliptic complex of differential operators
\be\label{inst-complex}
\begin{tikzcd}
\Lambda^0(\mathfrak{g}) \arrow{r}{d_A} & \Lambda^1(\mathfrak{g}) \arrow{r}{d^+_A} & \Lambda^2_+(\mathfrak{g})
\end{tikzcd}
\ee
Here $d_A$ is the exterior covariant derivative with respect to the background gauge field $A$, and the second operator is the restriction of the exterior covariant derivative (of a Lie algebra valued 1-form) to the space of self-dual Lie algebra valued 2-forms. The composition of the two operators appearing above vanishes when the background gauge field is an instanton $F(A)_+=0$. The above complex can be checked to be elliptic, which means that for the corresponding complex of symbols of the operators the image of the previous map coincides with the kernel of the next map. Also, the full YM equations $\star d_A \star F(A)=0$ are equivalent to $\star d_A \star F_+(A)=0$, and their linearisation is then the statement that 
\be\label{dd-star-YM}
(d^+_A)^* d^+_A a =0,
\ee
where $(d^+_A)^*$ is the adjoint of $d^+_A a$. So, the complex (\ref{inst-complex}) can be said to also know about the full second-order YM field equations. Adding the adjoint operators, one gets an elliptic operator
\be\label{dirac-YM}
\begin{tikzcd}
\Lambda^1(\mathfrak{g})\quad \arrow{r}{d^+_A+ (d_A)^*}  & \quad \Lambda^0(\mathfrak{g})  \oplus \Lambda^2_+(\mathfrak{g}).
\end{tikzcd}
\ee
This operator is known to be the (twisted by $A$) Dirac operator $D: S_+\otimes S_- \otimes \mathfrak{g} \to S_+\otimes S_+ \otimes \mathfrak{g}$, where $S_\pm$ are the two spinor representations of ${\rm Spin}(4)$. The instanton complex plays the key role in understanding the (local) properties of the moduli space of YM instantons, see \cite{Donaldson:1983wm}.

The aim of this paper is to build an analogous elliptic differential complex, but this time relevant for Einstein rather than YM equations. Our construction became possible due to the recent re-interpretation \cite{Bhoja:2024xbe} of the Pleba\'nski description \cite{Plebanski:1977zz} of Einstein equations. Let us remind the reader the basic points of this encoding of Einstein equations. 

Let $\Sigma^i, i=1,2,3\in \Lambda^2$ be a triple of 2-forms satisfying the condition $\Sigma^i\wedge \Sigma^j\sim \delta^{ij}$. Some literature refers to such a set of 2-forms as {\bf perfect triples}. The ${\rm GL}(4,\R)$ stabiliser of a perfect triple can be shown to be ${\rm SU}(2)$, and so perfect triples give an example of a $G$-structure, with $G={\rm SU}(2)$. Given an ${\rm SU}(2)$-structure $\Sigma^i$, there exists a triple of 1-forms $A^i$ satisfying the equation
\be\label{pleb-1}
d\Sigma^i + \epsilon^{ijk} A^j\wedge \Sigma^k=0.
\ee
This triple of 1-forms can be shown \cite{Bhoja:2024xbe} to coincide with the so-called intrinsic torsion of the ${\rm SU}(2)$-structure $\Sigma^i$. When the intrinsic torsion vanishes the ${\rm SU}(2)$-structure is integrable. It can be shown that in this case the metric defined by $\Sigma^i$ is Ricci-flat and there are 3 integrable complex structures satisfying the algebra of the imaginary quaternions. So, when $d\Sigma^i=0$ the metric defined by $\Sigma^i$ is hyper-K\"ahler. For more information about this description of hyper-K\"ahler metrics see e.g. \cite{Jiang}, Section 4.1. More generally, Einstein equations $Ric=0$ can also be encoded in the same language. Indeed, as is well-known from the general context of $G$-structures, the intrinsic torsion of a $G$-structure, when $G\subset{\rm SO}(n)$, where ${\rm dim}(M)=n$, together with its first derivatives, encode some part of the Riemann curvature tensor of the Riemannian metric defined by the $G$-structure. This is why it may be possible to impose some equations on the curvature, for example Einstein equations, as a first-order in derivatives condition on the intrinsic torsion. In the case at hand, one shows that Einstein equations $Ric=0$ become the statement that 
\be\label{pleb-2}
F^i(A) = \Psi^{ij} \Sigma^j, \qquad F^i(A) = dA^i + \frac{1}{2} \epsilon^{ijk} A^j\wedge A^k.
\ee
Here $\Psi^{ij}$ is an arbitrary symmetric tracefree $3\times 3$ matrix. 

The equations (\ref{pleb-1}), (\ref{pleb-2}) are the two equations of the Pleba\'nski formalism for General Relativity. Their linearisation at the background $d\Sigma^i=0$ (and thus $A^i=0$) is what gives rise to the Pleba\'nski complex that we are to define. Thus, we now consider a perfect triple of closed $d\Sigma^i=0$ 2-forms $\Sigma^i$. Define $E=\mathfrak{su}(2) = \mathfrak{so}(3)=\R^3$. Let us define $S\subset E\otimes \Lambda^2$ to be the tangent space to the space of perfect triples at $\Sigma^i$. We will characterise $S$ more explicitly below. Similarly, let $E\times \Lambda^1$ be the tangent space to the space of triples of 1-forms $A^i$. The linearisation of the equation (\ref{pleb-1}) at $A^i=0$ is then
\be\label{lin-pleb-1}
d \sigma^i + \epsilon^{ijk} a^j \wedge \Sigma^k=0, \qquad \sigma^i\in S, \quad a^i \in E\otimes \Lambda^1.
\ee
The solution of this equation can be stated in terms of a certain linear operator $J_1: E\otimes \Lambda^1 \to E\otimes \Lambda^1$ defined in (\ref{J-1}). This operator has two eigenspaces in $E\otimes \Lambda^1$, decomposing it $E\otimes \Lambda^1 = (E\otimes \Lambda^1)_4 \oplus (E\otimes \Lambda^1)_8$. The equation (\ref{lin-pleb-1}) becomes a linear equation for $a^i$, whose solution is given by $a^i = d_2 \sigma^i$, where 
\be\label{map-2}
d_2:S\to E\otimes \Lambda^1, \qquad d_2 \sigma^i = \frac{1}{2} J_1^{-1}( \star d \sigma^i).
\ee
Here $\star$ is the Hodge star with respect to the metric defined by $\Sigma^i$. This is the operator that appears as the second arrow of our to-be-constructed complex. 

To construct the first arrow we note that the objects $a^i$ that appear as the result of the map $d_2: S\to E\otimes \Lambda^1$ are diffeomorphism invariant. Diffeomorphisms act on the space $S$ via
\be\label{map-1}
 d_1: TM \to S, \qquad d_1 \xi = d (\xi\lrcorner  \Sigma^i).
\ee
 We will see in the main exposition that the object $d (\xi\lrcorner  \Sigma^i)$ is in $S$, the tangent space to the space of perfect triples. It is then immediate to see that the composition $d_2 d_1=0$. 
 
 The final arrow of the Pleba\'nski complex arises by considering the map 
 \be\label{map-3}
 d_3: E\otimes \Lambda^1 \to E, \qquad d_3 a^i = da^i \Big|_3 \equiv \epsilon^{ijk} \Sigma^j \wedge da^k / v_\Sigma.
 \ee
 Here $da^i\in E\otimes \Lambda^2$ is the exterior derivative of $a^i$, and the projection is taken onto the copy of $E$ inside $E\otimes \Lambda^2$. The last equality describes this projection explicitly, here $v_\Sigma$ is the volume form for the metric defined by $\Sigma^i$, which we can take to be $v_\Sigma = (1/6) \Sigma^i\wedge \Sigma^i$. The vanishing of the composition $d_3 d_2=0$ is the consequence of the equation (\ref{lin-pleb-1}). Indeed, taking the exterior derivative of this equation one gets $d_3 a^i=0$ whenever $a^i$ satisfies (\ref{lin-pleb-1}). 

 Assembling the spaces and the differential operators we get what we propose to call the Pleba\'nski complex
 \be\label{pleb-compl-intr}
 \begin{tikzcd}
TM  \arrow{r}{d_1} & S  \arrow{r}{d_2} & E\otimes \Lambda^1  \arrow{r}{d_3} & E
\end{tikzcd}
\ee
We then have the following
\begin{theorem} The complex (\ref{pleb-compl-intr}) is an elliptic complex of differential operators.
\end{theorem}
We remind the reader that a complex of first-order differential operators 
\be
 \begin{tikzcd}
 E_0  \arrow{r}{d_1} & E_1  \arrow{r}{d_2} & E_2  \arrow{r}{d_3} & E_3 
\end{tikzcd}
\ee
where $E_{0,1,2,3}$ are (the space of sections of) vector spaces, is said to be elliptic if the associated sequence of symbols
\be
 \begin{tikzcd}
0   \arrow{r}{} &  E_0  \arrow{r}{\sigma(d_1)} &  E_1  \arrow{r}{\sigma(d_2)} &  E_2  \arrow{r}{\sigma(d_3)} &  E_3  \arrow{r}{} & 0
\end{tikzcd}
\ee
is exact outside the zero section of $T^*M$. We remark that the dimensions of the spaces appearing in (\ref{pleb-compl-intr}) are
\be\label{pleb-compl-intr-dims}
 \begin{tikzcd}
4  \arrow{r}{d_1} & 13  \arrow{r}{d_2} & 12  \arrow{r}{d_3} & 3
\end{tikzcd}
\ee
and thus the alternating sum of the dimensions vanishes, as it should for an elliptic complex. 

Our second statement gives an Einstein analogue of the characterisation (\ref{dd-star-YM}) of the linearised YM equations.
\begin{theorem} Define the inner product on $E\otimes \Lambda^1$ to be given by
\be\label{inner-prod-a}
\langle a,a\rangle = \int_M \epsilon^{ijk} \Sigma^i \wedge a^j \wedge a^k.
\ee
The linearisation of the Einstein equations (\ref{pleb-2}) is the statement 
\be
d_2^* d_2 \sigma^i =0.
\ee
\end{theorem}
The operator $d_2^*$ also depends on the inner product in $S$, and this is not unique. But the above theorem holds for any choice of the inner product in $S$. 

Our final statement, whose proof will occupy the largest portion of the paper, is related to a construction of an elliptic operator for (\ref{pleb-compl-intr}) that squares to the Laplacian, giving a gravity analogue of (\ref{dirac-YM}). The most natural such construction would consist in simply adding the adjoints. Indeed, let us consider the operator
\be
D: S\oplus E \to TM \oplus E\otimes\Lambda^1, \qquad D(\sigma,\chi) = (d_1^*\sigma, d_2 \sigma + d_3^*\chi).
\ee
The operator so constructed is guaranteed to be an elliptic operator. However, somewhat to our surprise, we have the following negative result: $D^* D$ is {\bf not} a multiple of the Laplacian on $S\oplus E$, if we demand that the inner product on $E\otimes\Lambda^1$ is given by (\ref{inner-prod-a}). However, $D^* D$ is a multiple of the Laplacian if instead one uses the following inner product
\be
\langle a,a\rangle' = \int_M (a^i,a^i),
\ee
where $(\cdot,\cdot)$ is the metric pairing on $\Lambda^1$, defined with respect to the metric determine by $\Sigma^i$. While this choice of the inner product leads to the desired result for $D^* D$, this is not the choice that leads to $d_2^*$ being the operator appearing in the linearised Einstein equations. This makes us look for a construction of $D$ that has the property $D^* D\sim \Delta$ and that will involve the operators $d_2, d_2^*$, the later being defined with respect to the inner product (\ref{inner-prod-a}). Such a construction exists, but is non-trivial. 

Let us consider the more general elliptic differential operator 
\be
\tilde{D}: S\oplus E \to TM \oplus E\otimes\Lambda^1, \qquad   \tilde{D}(\sigma,\chi) = (\tilde{d}_1^* \sigma + \tilde{d}_4 \chi, d_2 \sigma + \tilde{d}_3^*\chi),
\ee
where we introduced $\tilde{d}_4: E\to TM$. We require that $d_2$ that appears here is the one of the Pleba\'nski complex (\ref{pleb-compl-intr}), and that the inner product on $E\otimes \Lambda^1$ is as in (\ref{inner-prod-a}) so that $d_2^*$ is the one that appears in the linearised Einstein equations. All other operators are allowed to be completely general first order differential operators as allowed by the relevant representation theory. In particular, we do not require that $\tilde{d}_1, d_2, \tilde{d}_3, \tilde{d}_4$ form a differential complex (i.e. there is no demand that the compositions are zero). We have the following
\begin{theorem} There is a choice of the inner product on $S$, and there is a choice of the operators $\tilde{d}_1, \tilde{d}_3, \tilde{d}_4$ such that all of the following hold:
\begin{itemize}
\item $\tilde{D}^* \tilde{D} =\pm \Delta$ on $S\oplus E$
\item The operators $\tilde{d}_1,\tilde{d}_3,\tilde{d}_4$ can be expressed in terms of $d_1, d_2, d_3$ of the Pleba\'nski complex. 
\item There exist linear maps
\be
T_1: S\oplus E\to S\oplus E, \qquad T_2 : TM\oplus E\otimes \Lambda^1 \to TM\oplus E\otimes \Lambda^1
\ee
such that 
\be\label{split}
T_2 \tilde{D} T_1 = D_4 \oplus D_{12},
\ee
with $D_4, D_{12}$ being operators
\be
D_4: \Lambda^0 \oplus E \to \Lambda^1, \qquad D_{12}: E \oplus {\rm Sym}^2_0(\Lambda^1) \to E\otimes \Lambda^1
\ee
explicitly given by (\ref{d4-d12}). Both of these operators are versions of the Dirac operator
\be
D_4: S_+\otimes S_+ \to S_-\otimes S_+, \qquad D_{12}: S_-\otimes S_-\otimes S_+^2 \to S_+\otimes S_-\otimes S_+^2.
\ee
\end{itemize}
\end{theorem}
Moreover, the choices we need to make for the above to hold are essentially unique, modulo signs can be absorbed into field redefinitions. The operator $\tilde{D}$ that is the result of construction in this theorem is explicitly given by
\be
\tilde{D} = \begin{pmatrix}\frac{1}{\sqrt{2}}(d^*_1 - \Phi d_2) && -\Phi d^*_3 \\ d_2 && \frac{1}{\sqrt{2}} J_1 d^*_3 \end{pmatrix}.
\ee
Here all operators are those of the Pleba\'nski complex, the linear map $\Phi: E\otimes \Lambda^1\to \Lambda^1$ is given by (\ref{Phi}) and $J_1: E\otimes \Lambda^1\to E\otimes \Lambda^1$ is the linear map given by (\ref{J-1}). This gives us a non-trivial gauge-fixing of the Pleba\'nski complex, by appropriately assembling the operators $d_1, d_2, d_3$ and their adjoints into an elliptic operator that can be written as (twisted by $T_1, T_2$) direct sum of two Dirac operators. This result requires quite a lot of involved calculations occupying Sections 6,7,8. 

If one does not demand that $\tilde{d}_1, \tilde{d}_3,\tilde{d}_4$ can be expressed in terms of $d_1, d_2, d_3$ of the Pleba\'nski complex then there is some more freedom in the choice of $\tilde{d}_1, \tilde{d}_3, \tilde{d}_4$ that leads to $\tilde{D}^* \tilde{D} =\pm \Delta$. The operators $D_4, D_{12}$ already appeared in the work \cite{Krasnov:2020bqr}, but this reference was using spinor techniques heavily, which makes it difficult to establish any uniqueness statement. In contrast, the present work makes all the assumptions that go into the construction of $\tilde{D}$ manifest, leading to the statement in the above theorem. 

We believe the operator $\tilde{D}$, together with demonstration that it splits as in (\ref{split}), to be the most interesting construction of this paper. In particular, we hope to develop a similar gauge-fixing for the full non-linear Einstein equations in a separate publication, having in mind applications in Numerical Relativity. The results of this paper are also likely to be useful for carrying out the Kummer construction of hyper-K\"ahler metrics on K3 surfaces, see the recent paper \cite{Jiang} for a review and new results in this direction. Another application of our construction is an analogous one in the context of holonomy ${\rm Spin}(7)$-manifolds, to be described elsewhere.

\section{\texorpdfstring{Decomposition of $E$-valued differential forms}{}}
\label{sec:decomp}

The material in this preparatory section is from \cite{Bhoja:2024xbe}. The rotation group in four dimensions is ${\rm SO}(4)={\rm SU}(2)\times{\rm SU}(2)/\Z_2$. One of the two ${\rm SU}(2)$ factors stabilises the 2-forms $\Sigma^i$. The other ${\rm SU}(2)$ factor rotates the triple $\Sigma^i$ by an ${\rm SO}(3)$ orthogonal transformation. This means that it is natural to treat $\Sigma^i$ as a map $\Sigma: \R^3 \to \Lambda^2$. The image of this map is the space $\Lambda^2_+$ of self-dual 2-forms. The map $\Sigma$ thus identifies $\R^3$ with $\Lambda^2_+$. This makes it natural to consider how the spaces of $E$-valued differential forms on $M$ transform with respect to all of the rotation group ${\rm SO}(4)={\rm SU}(2)\times{\rm SU}(2)/\Z_2$. 

Thus, we need to understand the decomposition of the spaces $\Lambda^1(M)\otimes E, \Lambda^2(M)\otimes E$, into irreducible representations of ${\rm SU}(2)\times{\rm SU}(2)/Z_2$. Irreducible representations of ${\rm SU}(2)$ are the spin $k/2$ representations that we denote by $S^k$. They are of dimension ${\rm dim}(S^k) = k+1$. As we have already discussed, there are two different ${\rm SU}(2)$'s in the game. One ${\rm SU}(2)$ is the group with respect to which the 2-forms $\Sigma^i$ are invariant. We will choose to denote this copy of ${\rm SU}(2)$ by ${\rm SU}_-(2)$, and the corresponding representations by $S^k_-$. The other ${\rm SU}(2)$ is one that acts non-trivially on $\Sigma^i$ by mixing them, with the map $\Sigma:E\to \Lambda^2$ being equivariant with respect to this copy of ${\rm SU}(2)$. We will denote it by ${\rm SU}_+(2)$, and the corresponding representations by $S_+^k$. We then have
\be
\Lambda^1(M) = S_+\otimes S_-, \qquad \Lambda^2(M) = S_+^2 \oplus S_-^2, \qquad E = S_+^2,
\ee
and the decomposition of $\Lambda^1(M)\otimes E, \Lambda^2(M)\otimes E$ into irreducibles is
\be\label{E-forms-decomp}
\Lambda^1(M) \otimes E= (S_+^3\otimes S_-) \oplus (S_+\otimes S_-), \\ \nonumber
\Lambda^2(M) \otimes E = S_+^4 \oplus S_+^2 \oplus C^\infty(M) \oplus (S_+^2 \otimes S_-^2).
\ee

\subsection{\texorpdfstring{Algebra of $\Sigma$'s}{}}

To obtain explicit formulas for the irreducible parts of $E$-valued differential forms, we need to define some linear operators acting on these spaces. These are built from the objects $\Sigma^i_{\mu\nu}$, together with the metric defined by $\Sigma$'s. We will be using index notation, and $\mu,\nu,\ldots=1,2,3,4$ are tensor indices.

The Riemannian metric defined by $\Sigma^i$ can be written explicitly as
\be\label{metric}
g_{\mu\nu} v_\Sigma = \frac{1}{6} \epsilon^{ijk} \Sigma^i_{\mu\alpha} \Sigma^j_{\nu\beta} \Sigma^k_{\gamma\delta} \tilde{\epsilon}^{\alpha\beta\gamma\delta}.
\ee
Here $\tilde{\epsilon}^{\alpha\beta\gamma\delta}$ is the densitiesed completely anti-symmetric tensor, which has components $\pm 1$ in any coordinate system and does not need a metric for its definition. Both sides of this formula are densitiesed $\mu\nu$ symmetric tensors, and $v_g=\sqrt{{\rm det}(g_{\mu\nu})}$ is the volume form for $g_{\mu\nu}$. 

One of the two indices of $\Sigma^i$ can be raised with the metric (\ref{metric}), to convert these objects into those in ${\rm End}(TM)$. We then have a triple of such endomorphisms of the tangent bundle, satisfying the algebra of the imaginary quaternions
\be\label{algebra}
\Sigma^i_{\mu}{}^\alpha \Sigma^j_\alpha{}^\nu = - \delta^{ij} \delta_\mu{}^\nu + \epsilon^{ijk} \Sigma^k_\mu{}^\nu.
\ee
There are also useful relations
\be\label{sigma-sigma}
\Sigma^i_{\mu\nu} \Sigma^i_{\rho\sigma}= g_{\mu\rho} g_{\nu\sigma} - g_{\mu\sigma} g_{\nu\rho} + \epsilon_{\mu\nu\rho\sigma}, \\ 
\label{sigma-sigma-epsilon}
\epsilon^{ijk} \Sigma^j_{\mu\nu} \Sigma^k_{\rho\sigma}= -2\Sigma^i_{[\mu|\rho|} g_{\nu]\sigma} + 2\Sigma^i_{[\mu|\sigma|} g_{\nu]\rho} , \\ \nonumber
\epsilon^{\mu\nu\rho\alpha} \Sigma^i_{\sigma\alpha} = \delta_\sigma^\rho \Sigma^{i\mu\nu} +  \delta_\sigma^\mu \Sigma^{i\nu\rho} + \delta_\sigma^\nu \Sigma^{i\rho\mu} .
\ee
These relations, as well as (\ref{algebra}), will be used on numerous occasions below, without explicit mention. 

\subsection{\texorpdfstring{Decomposition of $\Lambda^1(M) \otimes E$}{}}

Given an ${\rm SU}(2)$ structure $\Sigma^i$, we can define the following operator acting on $\Lambda^1(M) \otimes E$
\be\label{J-1}
\Lambda^1(M) \otimes E\ni a^i_\mu \to J_1(A)_\mu^i := \epsilon^{ijk} \Sigma^j_\mu{}^\alpha a^k_\alpha.
\ee
A simple calculation using (\ref{algebra}) shows that 
\be
J_1^2 = 2\mathbb{I} + J_1.
\ee
This means that the eigenvalues of $J_1$ are $2, -1$. The eigenspaces of $J_1$ are precisely the irreducibles appearing in the first line in (\ref{E-forms-decomp}). It is also easy to check that objects of the form
\be
\xi^\alpha \Sigma^i_{\alpha\mu} \in \Lambda^1(M)\otimes E
\ee
are eigenvectors of eigenvalue $2$. We then have the following characterisation 
\be\label{lambda-E-4}
(\Lambda^1(M)\otimes E)_4 = (S_+\otimes S_-) =\{ \xi^\alpha \Sigma^i_{\alpha\mu}, \xi\in TM\}.
\ee
The space 
\be\label{lambda-E-8}
(\Lambda^1(M)\otimes E)_8 = (S^3_+\otimes S_-) 
\ee
can then be characterised as the orthogonal complement of (\ref{lambda-E-4}) in $\Lambda^1(M)\otimes E$.

\subsection{\texorpdfstring{${\rm GL}(4)$ orbit of $\Sigma$ in $\Lambda^2(M) \otimes E$}{}}
To characterise some of the spaces appearing in the decomposition of $\Lambda^2(M) \otimes E$ we first consider the ${\rm GL}(4)$ orbit of the 2-forms $\Sigma^i$. The tangent space to this orbit is precisely the space $S$ of tangent vectors to perfect triples that we introduced above. This tangent space is the space of $E$-valued 2-forms of the form $h_{[\mu}{}^\alpha \Sigma^i_{ |\alpha|\nu]}$. Decomposing $h_{\mu\nu}\in{\rm GL}(4)$ into its symmetric and anti-symmetric parts, and noting that the anti-symmetric part is valued in $\Lambda^2(M)=S_+^2 \oplus S_-^2$, we get the following list of irreducibles appearing  
\be\label{S-space-param}
h_{[\mu}{}^\alpha \Sigma^i_{ |\alpha|\nu]} \in S= S_+^2 \oplus C^\infty(M) \oplus (S_+^2 \otimes S_-^2) \subset \Lambda^2(M) \otimes E,
\ee
which is all spaces in the second line of (\ref{E-forms-decomp}) apart from $S_+^4$. These irreducibles in $\Lambda^2(M) \otimes E$ can then be characterised as the images of the map $h_{\mu\nu}\to h_{[\mu}{}^\alpha \Sigma^i_{ |\alpha|\nu]} $. 

One can also act on the index $i$ of $\Sigma^i$ 2-forms with a ${\rm GL}(3)$ transformation, i.e., consider the orbit of $E$-valued 2-forms of the form $h^{ij} \Sigma^j_{\mu\nu}$. Decomposing the matrix $h^{ij}$ into symmetric and anti-symmetric parts, one finds the following list of irreducibles
\be
h^{ij} \Sigma^j_{\mu\nu} \in S_+^4 \oplus S_+^2 \oplus C^\infty(M).
\ee

In the opposite direction, given an object $B^i_{\mu\nu}\in \Lambda^2(M) \otimes E$, its irreducible parts can be extracted as follows
\be\label{irreducibles}
B^{(i}_{\alpha\beta} \Sigma^{j)\alpha\beta} - \frac{1}{3} \delta^{ij} B^{k}_{\alpha\beta} \Sigma^{k\alpha\beta} \in S_+^4, \\ \nonumber
\epsilon^{ijk} B^j_{\alpha\beta} \Sigma^{k\alpha\beta} \in S_+^2, \\ \nonumber
B^{k}_{\alpha\beta} \Sigma^{k\alpha\beta} \in C^\infty(M), \\ \nonumber
B^i_{\langle\mu|\alpha|} \Sigma^{i\alpha}{}_{\nu\rangle} \in S_+^2 \otimes S_-^2 .
\ee
Here $T_{\langle \mu\nu \rangle} = (1/2)(T_{\mu\nu} + T_{\nu\mu} - (1/4) g_{\mu\nu} g^{\alpha\beta} T_{\alpha\beta})$, which thus denotes the symmetric tracefree part. 

\subsection{Decomposition of \texorpdfstring{$\Lambda^2(M) \otimes E$}{}}
We can also describe the irreducible subspaces of $\Lambda^2(M) \otimes E$ as eigenspaces of a certain operator in $E$-valued 2-forms, similar to how we used $J_1$ to decompose $\Lambda^1\otimes E$. Let us introduce the following operator
\be\label{J2}
J_2: \Lambda^2\otimes E\to \Lambda^2\otimes E, \qquad J_2(B)_{\mu\nu}^i = \epsilon^{ijk} \Sigma^j_{[\mu}{}^\alpha B_{|\alpha|\nu]}^k, \qquad B_{\mu\nu}^i\in \Lambda^2\otimes E.
\ee
A computation gives
\be
J_2^2(B)_{\mu\nu}^i = \frac{1}{2} B_{\mu\nu}^i + \frac{1}{2}\epsilon_{\mu\nu}{}^{\alpha\beta} B^i_{\alpha\beta} + \frac{1}{2} J_2(B)_{\mu\nu}^i + \frac{1}{2} \Sigma^i_{[\mu}{}^\alpha \Sigma^j_{\nu]}{}^\beta B^j_{\alpha\beta}, \\ \nonumber
J_2^3(B)_{\mu\nu}^i =  \frac{1}{2}\epsilon_{\mu\nu}{}^{\alpha\beta} B^i_{\alpha\beta} + 2 J_2(B)_{\mu\nu}^i +  \Sigma^i_{[\mu}{}^\alpha \Sigma^j_{\nu]}{}^\beta B^j_{\alpha\beta}, \\ \nonumber
J_2^4(B)_{\mu\nu}^i =  \frac{1}{2} B_{\mu\nu}^i+ \frac{3}{2}\epsilon_{\mu\nu}{}^{\alpha\beta} B^i_{\alpha\beta} + \frac{5}{2} J_2(B)_{\mu\nu}^i +  \frac{5}{2} \Sigma^i_{[\mu}{}^\alpha \Sigma^j_{\nu]}{}^\beta B^j_{\alpha\beta}.
\ee
This implies
\be
J_2^4 - 2J_2^3-J_2^2 + 2 J_2=0 \qquad \text{or} \qquad J_2 (J_2-2)(J_2-1)(J_2+1)=0,
\ee
which implies that the eigenvalues of $J_2$ are $2,1,-1,0$. 

To characterise the eigenspaces we consider an arbitrary $3\times 3$ matrix $M^{ij}=M_s^{ij}+M_a^{ij}, M_s^{ij}=M_s^{(ij)}, M_a^{ij} = M_a^{[ij]}$ and compute
\be
J_2(M^{ij} \Sigma^j_{\mu\nu}) = {\rm Tr}(M) \Sigma^i_{\mu\nu}-M^{ji} \Sigma^j_{\mu\nu} = {\rm Tr}(M) \Sigma^i_{\mu\nu}-M_s^{ij} \Sigma^j_{\mu\nu}+ M_a^{ij} \Sigma^j_{\mu\nu}.
\ee
This means that the eigenspace of $J_2$ of eigenvalue $2$ is spanned by multiples of $\Sigma^i_{\mu\nu}$. The eigenspace of eigenvalue $1$ is $S_+^2$ spanned by $M_a^{ij} \Sigma^j_{\mu\nu}$. The eigenspace of eigenvalue $-1$ is $S_+^4$ spanned by $M_s^{ij} \Sigma^j_{\mu\nu}$ with ${\rm Tr}(M_s)=0$. 

We can also apply the operator $J_2$ to objects of the type $h_{[\mu}{}^\alpha \Sigma^i_{|\alpha|\nu]}$. We get
\be
J_2(h_{[\mu}{}^\alpha \Sigma^i_{|\alpha|\nu]}) = \frac{1}{2} h_\alpha{}^\alpha \Sigma^i_{\mu\nu}.
\ee
This means that the space $S_+^2\otimes S_-^2$ spanned by $h_{[\mu}{}^\alpha \Sigma^i_{|\alpha|\nu]}$ with tracefree $h_{\mu\nu}$ is eigenspace of $J_2$ of eigenvalue $0$.

All in all, we get
\be
\Lambda^2\otimes E= (\Lambda^2\otimes E)_5\oplus (\Lambda^2\otimes E)_3\oplus (\Lambda^2\otimes E)_1\oplus (\Lambda^2\otimes E)_9.
\ee
The last space here is $(\Lambda^2\otimes E)_9= \Lambda^-\otimes E$. The parametrisation of $S$ that we will use throughout the paper is
\be\label{param-sigma}
    \sigma^i_{\mu\nu} =  2 \epsilon^{ijk} \Sigma^j_{\mu\nu} h^k + 2 h_{[\mu}{}^\alpha \Sigma^i_{|\alpha|\nu]}
\ee
here $h_{\mu\nu}$ is a symmetric tracefree matrix and $h^i$ is an internal vector. Another notation that we will frequently use below is to represent objects in $S$ as the following two lists
\be
\left( h^i,\ h_{\mu\nu} \right) \in& E \oplus \rm Sym^2(\Lambda^1) \\
\left( h,\ h^i,\ \tilde{h}_{\mu\nu} \right) \in& C^\infty \oplus E \oplus \rm Sym^2_0(\Lambda^1)
\ee
where $\rm Sym^2_0(\Lambda^1)$ denotes the tracefree part and $h_{\mu\nu} = \frac{1}{4} g_{\mu\nu} h + \tilde{h}_{\mu\nu}$. So, the second representation splits the trace and the tracefree parts. 

\section{Proof of ellipticity}

We are now equipped to study the complex (\ref{pleb-compl-intr}) in more detail. We proof Theorem A by showing that the sequence of symbols 
\be\label{exact-sequence}
\begin{tikzcd}
0 \arrow{r}{} & TM  \arrow{r}{\sigma(d_1)} & S  \arrow{r}{\sigma(d_2)} & E\otimes \Lambda^1  \arrow{r}{\sigma(d_3)} & E \arrow{r}{} & 0
\end{tikzcd}
\ee
is exact, which means that Pleba\'nski complex is an elliptic complex of differential operators. 

\subsection{Proof of the theorem A}

We consider the symbols of all the differential operators. Symbols are obtained by replacing the operators of partial derivative $\partial_\mu$ with a factor of an arbitrary 1-form, which we denote by $k_\mu$.

The symbol of the operator in the first arrow (\ref{map-1}) is
\be
\xi^\mu \to k_{[\mu} \xi^\alpha \Sigma^i_{|\alpha|\nu]}.
\ee
To see that this has zero kernel we use the fact that if two 1-forms wedge to zero then they necessarily point in the same direction. The kernel of the above map is then 
\be\label{map-1-kernel-condition}
\xi^\alpha \Sigma^i_{\alpha\nu} = k_\nu \phi^i
\ee
where $\phi^i$ is an arbitrary internal vector. We would like to show that $\phi^i$ must vanish. 
We can solve this to find $\xi_\mu$ by contracting with $\Sigma^i_\mu{}^\nu$,
\be
\xi_\mu = \frac{1}{3} \Sigma^i_\mu{}^\nu \Sigma^i_{\alpha\nu} \xi^\alpha = \frac{1}{3} \Sigma^i_\mu{}^\nu k_\nu \phi^i.
\ee
Substituting this $\xi_\mu$ back into the left-hand side of (\ref{map-1-kernel-condition}) and using the algebra of $\Sigma^i$'s we find 
\be
\frac{1}{3} \Sigma^{j\alpha\nu} k_\nu \phi^j \Sigma^i_{\alpha\mu} = \frac{1}{3} k_\mu \phi^i + \frac{1}{3} \epsilon^{ijk} \Sigma^j_\mu{}^\nu k_\nu \phi^k.
\ee
Comparing with the right-hand side of (\ref{map-1-kernel-condition}) we find 
\be
2 k_\mu \phi^i = \epsilon^{ijk} \Sigma^j_\mu{}^\nu k_\nu \phi^k.
\ee
Contracting the above with $k^\mu$ we find 
\be
k^2 \phi^i = 0.
\ee
As $k^2 = k^\mu k_\mu \neq 0$ then this requires that $\phi^i = 0$.
Which implies that kernel of the first map (\ref{map-1}) is trivial, and 
shows the exactness of the first arrow in (\ref{exact-sequence}). 

The symbol of the second arrow (\ref{map-2}) is
\be\label{symbol-2}
\sigma^i \to \frac{1}{2} J_1^{-1} (\epsilon_\mu{}^{\alpha\beta\gamma} k_\alpha \sigma^i_{\beta\gamma}).
\ee
We know that $J_1^{-1}$ is invertible so we can safely ignore it and understand the kernel of 
\be
    \sigma^i \to \epsilon_\mu{}^{\alpha\beta\gamma} k_\alpha \sigma^i_{\beta\gamma}
\ee
Using the same linear algebra fact as before we know that $k_{[\alpha} \sigma^i_{\beta\gamma]}$ is zero when $\sigma^i_{\beta \gamma} = k_{[\beta} b^i_{\gamma]}$ for any $b^i_\mu \in E \otimes \Lambda^1$. This may make it seem that the kernel of the second map is all of $E \otimes \Lambda^1$. However, in general $k \wedge b^i \in E \times \Lambda^2$, which is bigger than $S$. Let us determine one-forms $b^i$ such that $k\wedge b^i \in S$. 
A 2-form from $E\otimes \Lambda^2$ belongs to $S$ if and only if 
\be\label{no-S-plus-4-in-kernel}
\Sigma^{<j|\mu\nu} \sigma^{i>}_{\mu\nu} = \Sigma^{<j|\mu\nu} k_\mu b^{|i>}_\nu = 0.
\ee
To find solutions of this equation, we decomposing $b^i$ into its irreducible parts
\be
b^i_\mu = \xi^\alpha \Sigma^i_{\alpha\mu} + (b_8)^i_\mu
\ee
where $\xi^\mu$ is the 4-vector part and $b^i_8\in (E\otimes\Lambda^1)_8$.
Using the algebra of $\Sigma^i$'s it is easy to see that the 4-vector part satisfies (\ref{no-S-plus-4-in-kernel}) automatically. Let us consider the $b^i_8$ part. To analyse this further we can choose a basis for $\Sigma^i$,
\be\label{self-dual-2-forms-k-basis}
\Sigma^i_{\mu\nu} = \frac{k_\mu}{|k|} e^i_\nu - \frac{k_\nu}{|k|} e^i_\mu - \epsilon^{ijk} e^j_\mu e^k_\nu
\ee
where $e^i_\mu$ are orthonormal and orthogonal to $k_\mu$, i.e. $e^i_\mu k^\mu = 0$ and $e^i_\mu e^{j\mu} = \delta^{ij}$.
This allows the decomposition 
\be
(b_8)^i_\mu = k_\mu b^i + b^{ij} e^j_\mu.
\ee
The eigenvalue of $b_8^i$ under the map $J_1$ is $-1$, which restricts the components of $b^{ij}$. Acting with $J_1$ gives
\be\label{J1-b}
J_1(b_8)^i_\mu =& \epsilon^{ijk}\left(\frac{k_\mu}{|k|} e^{j\nu} - \frac{k^\nu}{|k|} e^j_\mu - \epsilon^{jmn} e^m_\mu e^{n\nu}\right)(k_\nu b^k + b^{kl} e^l_\nu) \\ =& \frac{k_\mu}{|k|} \left( \epsilon^{ijk} b^{kj} \right) + \left( -b^{ji} + \delta^{ij} b^{kk} - |k| \epsilon^{ijk} b^k \right) e^j_\mu.
\ee
Therefore, $J_1(b_8)^i = -b^i_8$ implies
\be
 b^i =  \epsilon^{ijk} \frac{b^{jk}}{|k|} \quad \textrm{and} \quad b^{ij} = b^{ji} - \delta^{ij} b^{kk} + |k| \epsilon^{ijk} b^k.
\ee
The second equality implies $b^{ii}=0$, as well as the first. We have thus parametrised $b^i_8$ by a tracefree matrix $b^{ij}$ (not necessarily symmetric). 

We can now determine which part of $b^i_8$ is actually contributing to $k\wedge b^i_8\in S$. It is clear that the $k_\mu b^i$ part does not survive in $k\wedge b^i_8$. Let us consider the other part. We want to determine if any part of $k\wedge b^{ij} e^j$ is in $S$. Substituting this into (\ref{no-S-plus-4-in-kernel}) we find 
\be
\left(\frac{k^\mu}{|k|} e^{<i|\nu} - \frac{k^\nu}{|k|} e^{<i|\mu} - \epsilon^{<i|lm} e^{l\mu} e^{m\nu} \right) k_\mu b^{|j>n} e^n_\nu = |k| b^{<ij>}.
\ee
This must vanish in order for $k\wedge b^i_8\in S$. This means that $k\wedge b^i_8$ is either zero or not in $S$. It thus cannot contribute to the kernel of (\ref{symbol-2}). 
This means that $b^i_\mu = \xi^\alpha \Sigma^i_{\alpha\mu}$ generates all of the non-trivial kernel, where elements in the kernel are given by $\sigma^i_{\mu\nu} = k_{[\mu|} \xi^\alpha \Sigma^i_{\alpha |\nu]}$ and its dimension is $4$.
It is also clear then that the image of $\sigma(d_1)$ and the kernel of $\sigma(d_2)$ exactly coincide. 
As the dimension of the kernel of $\sigma(d_2)$ is $4$ and the dimension of the domain is $13$ then the image will have dimension $13-4 = 9$.

We now look at the final map
\be
a^i_\mu \to \epsilon^{ijk} \Sigma^{j\mu\nu} k_\mu a^k_\nu = k^\mu J_1(a)^i_\mu.
\ee
To find its kernel we decompose $a^i_\mu$
\be\label{E-valued-1-form-k-basis}
a^i_\mu = k_\mu a^i + a^{ij} e^j_\mu
\ee
where $a^i$ and $a^{ij}$ are an arbitrary internal vector and matrix.
Computing the action of $J_1$ is the same computation as (\ref{J1-b}), which gives
\be
J_1(a)^i_\mu =  \frac{k_\mu}{|k|}( - \epsilon^{ijk} a^{jk}) + \left( \delta^{ij} a^{kk} - a^{ji} - \epsilon^{ijk}|k| a^k \right) e^j_\mu.
\ee
This implies
\be
k^\mu J_1(a)^i_\mu = - |k| \epsilon^{ijk} a^{jk}.
\ee
This shows that the kernel of $\sigma(d_3)$ is all of $a^i_\mu$ apart from $a^{[ij]}$. So, the kernel is $12-3=9$ dimensional, and the image is $E$. 

We can also compute the image of $\sigma(d_2)$ explicitly, by using the form (\ref{d2-explicitly}) of this operator. In terms of symbols we have
\be\label{symbol-d2-explicitly}
(h^i , h_{\mu\nu}) \to - \Sigma^{i\alpha\beta} k_\alpha h_{\mu\beta} + 2 k_\mu h^i.
\ee
We now decompose $h_{\mu\nu}$ into its longitudinal and transverse parts
\be\label{h-decomp}
h_{\mu\nu} = \psi k_\mu k_\nu + k_{(\mu} X_{\nu)} + \tilde{h}_{\mu\nu}.
\ee
Here both $X_\mu, \tilde{h}_{\mu\nu}$ are transverse $k^\mu X_\mu=0, k^\mu \tilde{h}_{\mu\nu}=0$. We now want to compute the image of various part of the map (\ref{symbol-d2-explicitly}) and check that it coincides with the described above kernel of $\sigma(d_3)$. 

The last term in (\ref{symbol-d2-explicitly}) is precisely the first term in (\ref{E-valued-1-form-k-basis}). The image of the first term in (\ref{h-decomp}) is trivial. The image of the second term is
\be
- \frac{1}{2} \Sigma^{i\alpha\beta} k_\alpha X_\beta k_\mu,
\ee
which is also of the form of the first term in (\ref{E-valued-1-form-k-basis}). To compute the image of the last term in (\ref{h-decomp}) we use the explicit form of $\Sigma^{i\alpha\beta}$. We have
\be\label{image-d2-calc}
- \left(\frac{k^\alpha}{|k|} e^{i\beta} - \frac{k^\beta}{|k|} e^{i\alpha} - \epsilon^{ijk} e^{j\alpha} e^{k\beta}\right) k_\alpha \tilde{h}_{\mu\beta} = - |k| e^{i\beta} \tilde{h}_{\mu\beta}.
\ee
We can further decompose $\tilde{h}_{\mu\nu}$ into the basis of vectors $e^i_\mu$
\be
\tilde{h}_{\mu\nu} = h^{ij} e^i_\mu e^j_\nu,
\ee
where $h^{ij}=h^{(ij)}$ is a symmetric $3\times 3$ matrix. This gives for the right-hand-side of (\ref{image-d2-calc})
\be
- |k| h^{ij} e^j_\mu.
\ee
This is precisely the part of the second term in (\ref{E-valued-1-form-k-basis}) that is in the kernel of $\sigma(d_3)$. Thus, we have explicitly checked that the image of $\sigma(d_2)$ coincides with the kernel of $\sigma(d_3)$. This finishes the proof of ellipticity of the complex.

\section{Characterisation of Einstein equations}

The purpose of this section is to show that there exists a choice of the inner products on $E\otimes \Lambda^1$ and $S$ such that the linearisation of the Einstein equations $Ric=0$ becomes the statement $d_2^* d_2 \sigma^i=0$. 

\subsection{Linearisation of the Einstein condition}

As we have reviewed in the Introduction, in Pleba\'nski formalism the Einstein condition becomes the statement (\ref{pleb-2}). Its linearisation around $A^i=0$ is the statement
\be\label{einstein-lin}
d a^i = \psi^{ij} \Sigma^j,
\ee
with again $\psi^{ij}$ being an arbitrary symmetric tracefree $3\times 3$ matrix. From the discussion in subsection 2.4 we know that the right-hand side of (\ref{einstein-lin}) is an arbitrary vector in $(\Lambda^2\otimes E)_5\subset \Lambda^2\otimes E$. This means that (\ref{einstein-lin}) can be stated as $(da^i)\Big|_{1+3+9}=0$, i.e. as the statement that the projection of the exterior derivative of $a^i$ on the representations of dimensions $1,3,9$ in $\Lambda^2\otimes E$ vanishes. All these components can be recovered by computing
\be\label{einstein-lin-1}
(\partial_\mu a^i_\alpha - \partial_\alpha a^i_\mu) \Sigma^{i\alpha}{}_\nu \in \Lambda^1\otimes \Lambda^1.
\ee
The symmetric part of this tensor computes the $(\Lambda^2\otimes E)_{1+9}$ parts, and the anti-symmetric part only contains the $\Lambda^2_+ = (\Lambda^2\otimes E)_3$ part. So, linearised Einstein equations can be written as the statement that (\ref{einstein-lin-1}) vanishes. Projecting (\ref{einstein-lin-1}) onto $\Sigma^{k\mu\nu}$ it is easy to see that the $(\Lambda^2\otimes E)_3$ part of this condition reads
\be\label{einstein-lin-3}
\epsilon^{kij} \Sigma^{i\mu\nu} \partial_\mu a_\nu^j=0
\ee
The symmetric part of (\ref{einstein-lin-1}), on the other hand, gives
\be\label{einstein-lin-1+9}
\Sigma^i_{(\mu}{}^\alpha ( \partial_{\nu)} a^i_\alpha - \partial_\alpha a^i_{\nu)})=0.
\ee
As we already discussed, the condition (\ref{einstein-lin-3}) is an automatic consequence of the equation (\ref{lin-pleb-1}), and so holds automatically. The equations (\ref{einstein-lin-1+9}) are the ten linearised Einstein equations. 

\subsection{The adjoint of \texorpdfstring{$d_2$}{}}

We now compute the adjoint of $d_2$ given by (\ref{map-2}), with the inner product 
\be
\langle a,a\rangle = \int \epsilon^{ijk} \Sigma^{i\mu\nu} a_\mu^j a_\nu^k = - \int (a, J_1(a)),
\ee
where $J_1$ is the operator (\ref{J-1}) and $(a,b)=g^{\mu\nu} a_\mu^i b_\nu^i$ is the metric pairing of two objects in $\Lambda^1\otimes E$. We have
\be
\langle a, d_2 \sigma\rangle =& \langle a, \frac{1}{2} J_1^{-1} (\star d\sigma) \rangle= - \frac{1}{2} \int a_\mu^i \epsilon^{\mu\nu\rho\sigma} \partial_\nu \sigma^i_{\rho\sigma}
= -\frac{1}{2} \int \epsilon^{\mu\nu\rho\sigma} \sigma^i_{\mu\nu} \partial_\rho a_\sigma^i
\ee
Here we have used that $J_1^{-1} J_1=\mathbb{I}$. We already see that the adjoint of $d_2$ is related to the $da^i$. Given that $\sigma^i \in (\Lambda^2\otimes E)_{1+3+9}$, the operator $d_2^*$ will have components in all these spaces. It is clear that these will be multiples of the left-hand sides in (\ref{einstein-lin-3}) and (\ref{einstein-lin-1+9}). This proves Theorem B. We will present a more detailed computation of $d_2^*$ below, after the form of the inner product on $S$ is fixed. 

\section{The Pleba\'nski complex case}

We will now describe the Pleba\'nski complex in more details, including all the adjoint maps. We will be using the parametrisation (\ref{param-sigma}) of the space $S$, so that $S$ is the sum of the following irreducible components
\be\label{S-spaces}
S = \Lambda^0 \oplus E \oplus {\rm Sym}_0^2(\Lambda^1),
\ee
where ${\rm Sym}_0^2(\Lambda^1)$ denotes the tracefree part. We now establish the form of the operators $d_1, d_2, d_3$ in this parametrisation.

\subsection{Pleba\'nski complex in detail}

In the parametrisation (\ref{param-sigma}) of $S$ the first map of the Pleba\'nski complex takes the following form
\be
d_1 \xi = ( \partial^\mu \xi_\mu, \frac{1}{4} \Sigma^{i\mu\nu} \partial_\mu \xi_\nu, \frac{1}{2}( \partial_{\mu} \xi_{\nu}+ \partial_{\nu} \xi_{\mu}- \frac{1}{2} g_{\mu\nu} \partial^\alpha\xi_\alpha)).
\ee
Here the notation is that we list the irreducible components in $S$ in brackets, in the order they appear in (\ref{S-spaces}). 

The second operator is given by
\be
d_2 \sigma = \frac{1}{2} J_1^{-1}( \epsilon_\mu{}^{\alpha\beta\gamma} \partial_\alpha \sigma^i_{\beta\gamma}).
\ee
Let us compute this in the parametrisation (\ref{param-sigma}). We have
\be
\epsilon_\mu{}^{\alpha\beta\gamma} \partial_\alpha \sigma^i_{\beta\gamma} = 2 (\Sigma^{i\alpha\beta} \partial_\alpha h_{\mu\beta} - \Sigma^i_\mu{}^\beta \partial^\alpha h_{\alpha\beta} + \Sigma^i_\mu{}^\alpha \partial_\alpha h) + 4 \epsilon^{ijk} \Sigma^j_\mu{}^\alpha \partial_\alpha h^k.
\ee
Using $J_1^{-1} = (1/2)(J_1-\mathbb{I})$ we have
\be\nonumber
d_2 \sigma = \frac{1}{2} \epsilon^{ijk}\Sigma^j_\mu{}^\nu (\Sigma^{k\alpha\beta} \partial_\alpha h_{\nu\beta} - \Sigma^k_\nu{}^\beta \partial^\alpha h_{\alpha\beta} + \Sigma^k_\nu{}^\alpha \partial_\alpha h) - \frac{1}{2} (\Sigma^{i\alpha\beta} \partial_\alpha h_{\mu\beta} - \Sigma^i_\mu{}^\beta \partial^\alpha h_{\alpha\beta} + \Sigma^i_\mu{}^\alpha \partial_\alpha h)
\\ \nonumber 
+ \epsilon^{ijk}\Sigma^j_\mu{}^\nu \epsilon^{klm} \Sigma^l_\nu{}^\alpha \partial_\alpha h^m - \epsilon^{ijk} \Sigma^j_\mu{}^\alpha \partial_\alpha h^k.
\ee
The first term here is
\be
\frac{1}{2} \epsilon^{ijk}\Sigma^{j\mu\nu} \Sigma^{k\rho\sigma} \partial_\rho h_{\nu\sigma} = \frac{1}{2}  (- g^{\mu\rho} \Sigma^{i\nu\sigma} + g^{\nu\rho} \Sigma^{i\mu\sigma} + g^{\mu\sigma} \Sigma^{i\nu\rho} - g^{\nu\sigma} \Sigma^{i\mu\rho})  \partial_\rho h_{\nu\sigma}=\\ \nonumber
\frac{1}{2} \Sigma^{i\mu\sigma} \partial^\nu h_{\nu\sigma}+ \frac{1}{2} \Sigma^{i\nu\rho} \partial_\rho h_{\nu\mu} - \frac{1}{2} \Sigma^{i\mu\rho} \partial_\rho h.
\ee
The second term
\be
- \frac{1}{2} \epsilon^{ijk}\Sigma^j_\mu{}^\nu  \Sigma^k_\nu{}^\beta \partial^\alpha h_{\alpha\beta} =- \Sigma^i_\mu{}^\beta  \partial^\alpha h_{\alpha\beta} .
\ee
The third is
\be
\frac{1}{2} \epsilon^{ijk}\Sigma^j_\mu{}^\nu \Sigma^k_\nu{}^\alpha \partial_\alpha h =  \Sigma^i_\mu{}^\alpha  \partial_\alpha h.
\ee
The first term in the second line is
\be\nonumber
\epsilon^{ijk}\Sigma^j_\mu{}^\nu \epsilon^{klm} \Sigma^l_\nu{}^\alpha \partial_\alpha h^m = ( \delta^{il}\delta^{jm} - \delta^{im}\delta^{jl}) ( - \delta^{jl} \delta_\mu^\alpha + \epsilon^{jls}\Sigma^s_\mu{}^\alpha) \partial_\alpha h^m = 2 \partial_\mu h^i + \epsilon^{ijk} \Sigma^j_\mu{}^\alpha \partial_\alpha h^k.
\ee
Combining everything we observe multiple cancellations, with the final result being
\be\label{d2-explicitly}
d_2 \sigma = -  \Sigma^{i\alpha\beta} \partial_\alpha h_{\mu\beta} +2 \partial_\mu h^i.
\ee
Let us now exhibit the tracefree part and the trace of $h_{\mu\nu}$ explicitly. We parametrise
\be\label{h-tilde}
h_{\mu\nu} = \tilde{h}_{\mu\nu} + \frac{1}{4} g_{\mu\nu} h.
\ee
In this parametrisation 
\be
d_2 \sigma = \frac{1}{4} \Sigma^i_\mu{}^\alpha \partial_\alpha h +2 \partial_\mu h^i -  \Sigma^{i\alpha\beta} \partial_\alpha \tilde{h}_{\mu\beta} .
\ee

The last operator of the complex is
\be
d_3 a = \epsilon^{ijk} \Sigma^{j\mu\nu} \partial_\mu a_\nu^k.
\ee

\subsection{Inner products}

We will now write the general form of the inner products in all the spaces. We will be identifying $TM\sim \Lambda^1$ using the metric. The inner products on $\Lambda^1$ and $E$ are unique (up to an overall multiple which we choose to be unity). The spaces $S, E\otimes \Lambda^1$ each contains several irreducible components, and the inner product is no longer unique. We will introduce arbitrary (at this stage) constants to parametrise the arising inner products.

The inner product in $\Lambda^1$ is
\be
\langle \xi,\xi\rangle = \int_M (\xi_\mu)^2
\ee
where $(\xi_\mu)^2=\xi^\mu \xi_\mu$ denotes the metric contraction. The inner product in $S$ an arbitrary linear combination of inner products on all 3 irreducible pieces composing $S$
\be
\langle \sigma,\sigma \rangle = \int_M \beta_1 h^2 + \beta_2 (h^i)^2 + \beta_3 (\tilde{h}_{\mu\nu})^2.
\ee
The inner product in $E\otimes \Lambda^1$ is similarly a linear combination of two terms
\be
\langle a,a\rangle = \int_M \gamma_1 (a_\mu^i)^2 + \gamma_2 \epsilon^{ijk} \Sigma^{i\mu\nu} a^j_\mu a^k_\nu.
\ee
The inner product in $E$ is
\be
\langle \chi,\chi\rangle = \int_M (\chi^i)^2.
\ee

\subsection{General form of the adjoints}

We will first write down the most general form that the adjoint operators can take (based on the relevant representation theory). For the adjoint of the first operator we have
\be\label{d1-adj-gen}
d_1^* \sigma = a_1' \partial_\mu h + a_2' \Sigma^i_\mu{}^\nu \partial_\nu h^i + a_3' \partial^\nu \tilde{h}_{\mu\nu}.
\ee
For the second operator we have
\be\label{d2-adj-gen}
d_2^* a = (b_1' \Sigma^{i\mu\nu} \partial_\mu a_\nu^i,\ b_2' \partial^\mu a_\mu^i + b_3' \epsilon^{ijk} \Sigma^{j\mu\nu} \partial_\mu a_\nu^k,\ 
2b_4' \Sigma^i_{<\mu}{}^\alpha \partial_{\nu>}a_\alpha^i +  2b_5' \Sigma^i_{<\mu}{}^\alpha \partial_\alpha a_{\nu>}^i).
\ee
The third adjoint has the form
\be\label{d3-adj-gen}
d_3^* \chi = c_1' \partial_\mu \chi^i + c_2' \epsilon^{ijk} \Sigma^j_\mu{}^\alpha \partial_\alpha \chi^k.
\ee

\subsection{Pleba\'nski case adjoints}
We now compute the Pleba\'nski complex adjoints. For $\langle \sigma, d_1 \xi\rangle = \langle d_1^* \sigma, \xi\rangle$ we get
\be
\int_M \beta_1 h \partial^\mu \xi_\mu + \frac{\beta_2}{4}  h^i \Sigma^{i\mu\nu} \partial_\mu \xi_\nu + \beta_3  \tilde{h}^{\mu\nu} \partial_\mu \xi_\nu =\\ \nonumber
\int_M \xi^\mu (a_1' \partial_\mu h + a_2' \Sigma^i_\mu{}^\nu \partial_\nu h^i + a_3' \partial^\nu \tilde{h}_{\mu\nu}).
\ee
So, we have
\be
a_1'= -\beta_1 , \qquad  a_2'= \frac{\beta_2}{4} , \qquad a_3'= - \beta_3.
\ee

For $\langle a, d_2 \sigma \rangle = \langle d_2^* a, \sigma \rangle$ we have
\be\nonumber
\int_M \gamma_1 a^{i\mu} (  \frac{1}{4} \Sigma^i_\mu{}^\alpha \partial_\alpha h +2 \partial_\mu h^i -  \Sigma^{i\alpha\beta} \partial_\alpha \tilde{h}_{\mu\beta} )
+ \gamma_2 \epsilon^{ijk} \Sigma^{i\mu\nu} a^j_\mu ( \frac{1}{4} \Sigma^k_\nu{}^\alpha \partial_\alpha h +2 \partial_\nu h^k -  \Sigma^{k\alpha\beta} \partial_\alpha \tilde{h}_{\nu\beta}  ) = \\ \nonumber
\int_M \beta_1 b_1' h \Sigma^{i\mu\nu} \partial_\mu a_\nu^i + \beta_2 h^i (b_2' \partial^\mu a_\mu^i + b_3' \epsilon^{ijk} \Sigma^{j\mu\nu} \partial_\mu a_\nu^k) + \beta_3 \tilde{h}^{\mu\nu} (2b_4' \Sigma^i_{(\mu}{}^\alpha \partial_{\nu)}a_\alpha^i +  2b_5' \Sigma^i_{(\mu}{}^\alpha \partial_\alpha a_{\nu)}^i).
\ee
The last term in the first line gives
\be
\epsilon^{ijk} \Sigma^{i\mu\nu}   \Sigma^{k\rho\sigma} a^j_\mu \partial_\rho \tilde{h}_{\nu\sigma} = ( g^{\mu\rho} \Sigma^{i\nu\sigma} - g^{\nu\rho} \Sigma^{i\mu\sigma} - g^{\mu\sigma} \Sigma^{i\nu\rho} + g^{\nu\sigma} \Sigma^{i\mu\rho})a^i_\mu \partial_\rho \tilde{h}_{\nu\sigma} = \\ \nonumber
- \Sigma^{i\mu\sigma}a^i_\mu \partial^\nu \tilde{h}_{\nu\sigma} - \Sigma^{i\nu\rho} a^{i\mu} \partial_\rho \tilde{h}_{\nu\mu},
\ee
where we used the fact that $\tilde{h}_{\mu\nu}$ is traceless. This means we have
\be\nonumber
\int_M \frac{\gamma_1- 2\gamma_2}{4} \Sigma^{i\mu\rho}a^i_\mu \partial_\rho h +2\gamma_1 a^{i\mu} \partial_\mu h^i + 2\gamma_2 \epsilon^{ijk} \Sigma^{i\mu\nu} a^j_\mu  \partial_\nu h^k 
- (\gamma_1 +\gamma_2)a^{i\mu} \Sigma^{i\alpha\beta} \partial_\alpha \tilde{h}_{\mu\beta} +\gamma_2 \Sigma^{i\mu\sigma}a^i_\mu \partial^\nu \tilde{h}_{\nu\sigma} 
 = \\ \nonumber
\int_M \beta_1 b_1' h \Sigma^{i\mu\nu} \partial_\mu a_\nu^i + \beta_2 h^i (b_2' \partial^\mu a_\mu^i + b_3' \epsilon^{ijk} \Sigma^{j\mu\nu} \partial_\mu a_\nu^k) + \beta_3 \tilde{h}^{\mu\nu} (2b_4' \Sigma^i_{(\mu}{}^\alpha \partial_{\nu)}a_\alpha^i +  2b_5' \Sigma^i_{(\mu}{}^\alpha \partial_\alpha a_{\nu)}^i).
\ee
This gives
\be
b_1' =  \frac{\gamma_1-2\gamma_2}{4\beta_1}, \quad b_2' = - \frac{2\gamma_1}{\beta_2}, \quad b_3'= \frac{2\gamma_2}{\beta_2}, \quad b_4'= \frac{\gamma_2}{2\beta_3}, \quad b_5'= - \frac{\gamma_1+\gamma_2}{2\beta_3}.
\ee

For $\langle \chi, d_3 a\rangle = \langle d_3^* \chi, a\rangle$ we have
\be
\int_M  \chi^i  \epsilon^{ijk} \Sigma^{j\mu\nu} \partial_\mu a_\nu^k  =\\ \nonumber \int_M 
\gamma_1 a^{i\mu} (c_1' \partial_\mu \chi^i + c_2' \epsilon^{ijk} \Sigma^j_\mu{}^\alpha \partial_\alpha \chi^k) + \gamma_2 \epsilon^{ijk} \Sigma^{i\mu\nu} a_\mu^j (c_1' \partial_\nu \chi^k + c_2' \epsilon^{klm} \Sigma^l_\nu{}^\alpha \partial_\alpha \chi^m), 
\ee
or 
\be
\int_M  \chi^i  \epsilon^{ijk} \Sigma^{j\mu\nu} \partial_\mu a_\nu^k  =\\ \nonumber \int_M 
(\gamma_1 c_1' - 2 \gamma_2 c_2' ) a^{i\mu}  \partial_\mu \chi^i + (\gamma_1 c_2'  - \gamma_2 c_1' - \gamma_2 c_2' ) \epsilon^{ijk} a^{i\mu} \Sigma^j_\mu{}^\alpha \partial_\alpha \chi^k ,
\ee
and so
\be
\gamma_1 c_1' - 2 \gamma_2 c_2' =0, \qquad 1  = -\gamma_1 c_2'  + \gamma_2 c_1' + \gamma_2 c_2' .
\ee
This gives
\be
c_1' = \frac{2\gamma_2}{\gamma_1(\gamma_2-\gamma_1) + 2\gamma_2^2}, \quad c_2' = \frac{\gamma_1}{\gamma_1(\gamma_2-\gamma_1) + 2\gamma_2^2}.
\ee
We have explicitly checked that the compositions of the adjoint operators (with these coefficients) vanish, as they should. 

 \subsection{Imposing the \texorpdfstring{$\Delta$}{} condition}
 
 We now form the operators $D$ and $D^*$, see (\ref{D}) and (\ref{D*}). Requiring $D^*D$ to be a multiple of the $\Delta=\partial^\mu\partial_\mu$ operator, we are led to equations (\ref{D2-eqs}). Substituting the values of all the coefficients as they are (\ref{Pleb-coeffs}) in the Pleba\'nski complex case, and expressing the coefficients of the adjoint operators via those in the inner products, we can view the equations  (\ref{D2-eqs}) as those on the inner product coefficients. The solution is as follows
 \be
\langle \sigma,\sigma \rangle = \beta \int_M \frac{1}{4} h^2 + 8 (h^i)^2 +  (\tilde{h}_{\mu\nu})^2, \qquad
\langle a,a\rangle = \beta^2 \int_M  (a_\mu^i)^2 .
\ee
It also makes sense to choose $\beta=1$. We note that the relative coefficient in the $h^2, (\tilde{h}_{\mu\nu})^2$ terms is the natural one. We can pass to treating $h_{\mu\nu}$ as an object that contains the trace part. Indeed, using (\ref{h-tilde}) we have
\be\label{trace-tracefree}
(h_{\mu\nu})^2 = (\tilde{h}_{\mu\nu})^2 + \frac{1}{4} h^2, 
\ee
which is precisely what we have in the above inner product. So, we can write the arising via this computation inner products as
\be\label{inner-prod-1}
\langle \xi,\xi \rangle =\int_M (\xi_\mu)^2, \quad
\langle \sigma,\sigma \rangle = \int_M 8 (h^i)^2 +  (h_{\mu\nu})^2, \quad
\langle a,a\rangle =  \int_M  (a_\mu^i)^2 , \quad
\langle \xi,\xi\rangle =  \int_M  (\xi^i)^2.
\ee

\subsection{Inner product on \texorpdfstring{$S$ in terms of $\sigma^i_{\mu\nu}$}{}} 

It is interesting to rewrite the inner product on $S$ in (\ref{inner-prod-1}) in terms of the object $\sigma^i_{\mu\nu}\in S$. Using the parametrisation (\ref{param-sigma}) we have
\be\nonumber
(\sigma^i_{\mu\nu})^2 = 2 (h_\mu{}^\alpha \Sigma^i_{\alpha\nu} - h_\nu{}^\alpha \Sigma^i_{\alpha\mu}) h^{\mu\beta} \Sigma^i_\beta{}^\nu + 8 h_\mu{}^\alpha \Sigma^i_{\alpha\nu} \epsilon^{ijk} \Sigma^{j\mu\nu} \chi^k + 4 \epsilon^{ijk} \Sigma^j_{\mu\nu} \chi^k \epsilon^{imn} \Sigma^{m\mu\nu} \chi^n= \\ \nonumber
2( 3 (h_{\mu\nu})^2 - h_\nu{}^\alpha h^{\mu\beta} ( g_{\alpha\beta} \delta_\mu{}^\nu - \delta_\alpha{}^\nu g_{\mu\beta} + \epsilon_{\alpha\mu\beta}{}^\nu)) + 32 (\chi^i)^2 = \\ \nonumber
2( 2(h_{\mu\nu})^2 + h^2) + 32 (\chi^i)^2.
\ee
In the parametrisation (\ref{h-tilde}) this becomes
\be
\frac{1}{4} (\sigma^i_{\mu\nu})^2 = (\tilde{h}_{\mu\nu})^2 + \frac{3}{4} h^2 + 8 (\chi^i)^2.
\ee
We also have
\be
\Sigma^{i\mu\nu} \sigma^i_{\mu\nu} = 2 \Sigma^{i\mu\nu} h_\mu{}^\alpha \Sigma^i_{\alpha\nu} = 6h.
\ee
This means that 
\be
\frac{1}{4} h^2 + 8 (h^i)^2 +  (\tilde{h}_{\mu\nu})^2 = \frac{1}{4} (\sigma^i_{\mu\nu})^2 - \frac{1}{72} ( \Sigma^{i\mu\nu} \sigma^i_{\mu\nu} )^2.
\ee

\subsection{Explicit form of the adjoints}

We now write the Pleba\'nski complex operators without splitting the tracefree and the trace parts of $h_{\mu\nu}$
\be
d_1 \xi = ( \frac{1}{4} \Sigma^{i\mu\nu} \partial_\mu \xi_\nu, \partial_{(\mu} \xi_{\nu)}) \in (E,S^2 T^*M) , \\ \nonumber
d_2 \sigma = -  \Sigma^{i\alpha\beta} \partial_\alpha h_{\mu\beta} +2 \partial_\mu h^i, \\ \nonumber
d_3 a = \epsilon^{ijk} \Sigma^{j\mu\nu} \partial_\mu a_\nu^k.
\ee
The adjoints, in the same notation are given by 
\be
d_1^* \sigma = 2 \Sigma^i_\mu{}^\nu \partial_\nu h^i - \partial^\nu h_{\mu\nu}, \\ \nonumber
d_2^* a =( -\frac{1}{4} \partial^\mu a_\mu^i , - \Sigma^i_{(\mu}{}^\alpha \partial_\alpha a_{\nu)}^i) \in (E,S^2 T^*M), \\ \nonumber
d_3^* \chi =  - \epsilon^{ijk} \Sigma^j_\mu{}^\alpha \partial_\alpha \chi^k.
\ee
As a check, we recompute the operators that appear in $D^* D$. The operator $(d_1 d_1^* +d_2^* d_2 ) \sigma$ becomes, in the $h^i$ component
\be
 =  \frac{1}{4} \Sigma^{i\mu\nu} \partial_\mu (2 \Sigma^j_\nu{}^\alpha \partial_\alpha h^j - \partial^\alpha h_{\nu\alpha})
 -\frac{1}{4} \partial^\mu (-  \Sigma^{i\alpha\beta} \partial_\alpha h_{\mu\beta} +2 \partial_\mu h^i) = \\ \nonumber
 -\frac{1}{2}   \partial^\mu \partial_\mu h^i -  \frac{1}{4} \Sigma^{i\mu\nu} \partial_\mu\partial^\alpha h_{\nu\alpha}
 + \frac{1}{4}   \Sigma^{i\alpha\beta} \partial^\mu \partial_\alpha h_{\mu\beta} - \frac{1}{2} \partial^\mu\partial_\mu h^i = -  \partial^\mu\partial_\mu h^i .
 \ee
 In the $h_{\mu\nu}$ component we get
 \be
 \partial_{(\mu} (2 \Sigma^i_{\nu)}{}^\alpha \partial_\alpha h^i - \partial^\alpha h_{\nu)\alpha})- \Sigma^i_{(\mu}{}^\alpha \partial_\alpha  (-  \Sigma^{i\rho\sigma} \partial_\rho h_{\nu)\sigma} +2 \partial_{\nu)} h^i)= \\ \nonumber
 2 \Sigma^i_{(\mu}{}^\alpha \partial_{\nu)} \partial_\alpha h^i -  \partial_{(\mu} \partial^\alpha h_{\nu)\alpha} + ( \delta_{(\mu}{}^\rho g^{\alpha\sigma} - \delta_{(\mu}{}^\sigma g^{\alpha\rho}) \partial_\alpha \partial_\rho h_{\nu)\sigma} - 2\Sigma^i_{(\mu}{}^\alpha \partial_\alpha \partial_{\nu)} h^i = - \partial^\alpha \partial_\alpha h_{\mu\nu}.
 \ee
 For the operator $d_3 d_3^* \chi$ we have
 \be
 d_3 d_3^* \chi= \epsilon^{ijk} \Sigma^{j\mu\nu} \partial_\mu (- \epsilon^{klm} \Sigma^l_\nu{}^\alpha \partial_\alpha \chi^m)=
 - \epsilon^{ijk}\epsilon^{klm} (- \delta^{jl} g^{\mu\alpha}) \partial_\mu\partial_\alpha \chi^m = - 2 \partial^\mu \partial_\mu \chi^i,
 \ee
 which confirms that $D^* D\sim \Delta$.  
 
\subsection{Linearised Einstein equations}

As we have seen previously, see (\ref{einstein-lin-1+9}), the linearisation of the Einstein equations is 
\be\label{correct-d2-star}
\Sigma^i_{(\mu}{}^\alpha (\partial_{\nu)} a_\alpha^i - \partial_\alpha a_{\nu)}^i) =0.
\ee
We would like to realise the operator that appears in this equation as $d_2^*$, but this is not the operator $d_2^*$ that arises if the inner products are chosen as in (\ref{inner-prod-1}) to satisfy $D^*D\sim \Delta$ condition. We now change the inner product to obtain the desired $d_2^*$, but this will then be in conflict with the $D^*D\sim \Delta$ condition. Still, we develop this option as the one that is relevant for the Einstein equations. The construction of an elliptic operator $\tilde{D}$ with $\tilde{D}^*\tilde{D}\sim \Delta$ will occupy us in the next sections. 

\subsection{Pleba\'nski case different inner product}

We now consider the case of Pleba\'nski complex with the inner product $\gamma_1=0$. In this case $b_3=b_4=0$, and it is not hard to see that this gives $b_4'+b_5'=0, b_2'=0$. The other coefficients are
\be
b_1'= -\frac{\gamma_2}{2\beta_1}, \quad b_3'= \frac{2\gamma_2}{\beta_2}, \quad b_4'=-b_5'= \frac{\gamma_2}{2\beta_3}
\ee
so that 
\be\nonumber
d_2^* a = \gamma_2 (-\frac{1}{2\beta_1} \Sigma^{i\mu\nu} \partial_\mu a_\nu^i,  \frac{2}{\beta_2} \epsilon^{ijk} \Sigma^{j\mu\nu} \partial_\mu a_\nu^k, \frac{1}{\beta_3}(\Sigma^i_{(\mu}{}^\alpha \partial_{\nu)}a_\alpha^i -\Sigma^i_{(\mu}{}^\alpha \partial_\alpha a_{\nu)}^i) - \frac{1}{2\beta_3} g_{\mu\nu} \Sigma^{i\rho\sigma} \partial_\rho a_\sigma^i).
\ee
Note that the desired operator (\ref{correct-d2-star}) appears here.

We now choose 
\be
\gamma_2=1, \qquad \beta_1=\frac{1}{4}, \qquad \beta_2 = 8, \qquad \beta_3=1,
\ee
which corresponds to the following inner products
\be\label{inner-prod-pleb}
\langle \xi,\xi \rangle =\int_M (\xi_\mu)^2, \quad
\langle \sigma,\sigma \rangle = \int_M \frac{1}{4} h^2 +8 (h^i)^2 +  (\tilde{h}_{\mu\nu})^2, \\ \nonumber
\langle a,a\rangle =  \int_M  \epsilon^{ijk} \Sigma^{i\mu\nu} a_\mu^j a_\nu^k, \quad
\langle \xi,\xi\rangle =  \int_M  (\xi^i)^2.
\ee
These are the same inner products in all the spaces as in (\ref{inner-prod-1}), apart from that in $E\otimes\Lambda^1$, where we instead take the inner product that contains the $J_1$ operator. With these inner products, the full complex becomes
\be
d_1 \xi = ( \partial^\mu \xi_\mu, \frac{1}{4} \Sigma^{i\mu\nu} \partial_\mu \xi_\nu, \partial_{\langle\mu} \xi_{\nu\rangle}) \in (E,S^2 T^*M) , \\ \nonumber
d_2 \sigma = \frac{1}{4} \Sigma^i_\mu{}^\nu \partial_\nu h +2 \partial_\mu h^i -  \Sigma^{i\alpha\beta} \partial_\alpha \tilde{h}_{\mu\beta} , \\ \nonumber
d_3 a = \epsilon^{ijk} \Sigma^{j\mu\nu} \partial_\mu a_\nu^k,
\ee
and
\be\nonumber
d_1^* \sigma =   -\frac{1}{4}\partial_\mu h + 2\Sigma^i_\mu{}^\nu \partial_\nu h^i -  \partial^\nu \tilde{h}_{\mu\nu}, \\ \nonumber
d_2^* a = ( -2 \Sigma^{i\mu\nu} \partial_\mu a_\nu^i, \frac{1}{4} \epsilon^{ijk} \Sigma^{j\mu\nu} \partial_\mu a_\nu^k, \Sigma^i_{\langle\mu}{}^\alpha \partial_{\nu\rangle}a_\alpha^i-\Sigma^i_{\langle\mu}{}^\alpha \partial_\alpha a_{\nu\rangle}^i), \\ \nonumber
d_3^* \chi =   \partial_\mu \chi^i.
\ee
We will refer to the operators appearing above as those of the Pleba\'nski complex. Here $\langle\rangle$ denotes the symmetric tracefree part. 

\subsection{Computation of \texorpdfstring{$d_2^* d_2$}{}}

As a check, we compute the operator $d_2^* d_2$ for the operators of the Pleba\'nski complex as detailed above. We will need the result of this computation below. For the trace part we have
\be
-2 \Sigma^{i\mu\nu} \partial_\mu (\frac{1}{4} \Sigma^i_\nu{}^\alpha \partial_\alpha h +2 \partial_\nu h^i -  \Sigma^{i\rho\sigma} \partial_\rho \tilde{h}_{\nu\sigma} )=
\frac{3}{2} \partial^\mu \partial_\mu h - 2 \partial^\mu \partial^\nu \tilde{h}_{\mu\nu}.
\ee
For the vector part we have
\be
\frac{1}{4} \epsilon^{ijk} \Sigma^{j\mu\nu} \partial_\mu (\frac{1}{4} \Sigma^k_\nu{}^\alpha \partial_\alpha h +2 \partial_\nu h^k -  \Sigma^{k\rho\sigma} \partial_\rho \tilde{h}_{\nu\sigma} )=
\\ \nonumber
- \frac{1}{4}( g^{\mu\rho} \Sigma^{i\nu\sigma} - g^{\nu\rho} \Sigma^{i\mu\sigma} - g^{\mu\sigma} \Sigma^{i\nu\rho} + g^{\nu\sigma} \Sigma^{i\mu\rho}) \partial_\mu\partial_\rho \tilde{h}_{\nu\sigma} =0.
\ee
For the tracefree part we have
\be\nonumber
\Sigma^i_{\langle\mu}{}^\alpha \partial_{\nu\rangle}(\frac{1}{4} \Sigma^i_\alpha{}^\beta \partial_\beta h +2 \partial_\alpha h^i -  \Sigma^{i\rho\sigma} \partial_\rho \tilde{h}_{\alpha\sigma} )
-\Sigma^i_{\langle\mu}{}^\alpha \partial_\alpha (\frac{1}{4} \Sigma^i_{\nu\rangle}{}^\beta \partial_\beta h +2 \partial_{\nu\rangle} h^i -  \Sigma^{i\rho\sigma} \partial_\rho \tilde{h}_{\nu\rangle\sigma} )= \\ \nonumber
- \frac{3}{4} \partial_{\langle\mu} \partial_{\nu\rangle} h + 2 \Sigma^i_{\langle\mu}{}^\alpha \partial_{\nu\rangle}\partial_\alpha h^i + \partial_{\langle\mu}\partial^\rho \tilde{h}_{\nu\rangle \alpha}
+ \frac{1}{4} \partial_{\langle\mu} \partial_{\nu\rangle} h - 2 \Sigma^i_{\langle\mu}{}^\alpha \partial_{\nu\rangle}\partial_\alpha h^i + \partial_{\langle\mu}\partial^\rho \tilde{h}_{\nu\rangle \alpha}- \partial^\alpha \partial_\alpha \tilde{h}_{\mu\nu}= \\ \nonumber
-  \frac{1}{2} \partial_{\langle\mu} \partial_{\nu\rangle} h+ 2\partial_{\langle\mu}\partial^\rho \tilde{h}_{\nu\rangle \alpha}- \partial^\alpha \partial_\alpha \tilde{h}_{\mu\nu}.
\ee
Collecting the results we have
\be\label{d2}
d_2^* d_2 \sigma = ( \frac{3}{2} \partial^\mu \partial_\mu h - 2 \partial^\mu \partial^\nu \tilde{h}_{\mu\nu}, 0, -  \frac{1}{2} \partial_{\langle\mu} \partial_{\nu\rangle} h+ 2\partial_{\langle\mu}\partial^\rho \tilde{h}_{\nu\rangle \alpha}- \partial^\alpha \partial_\alpha \tilde{h}_{\mu\nu}).
\ee

\section{General complex}

We now change the game. We start by defining a general set of differential operators between the same spaces as appear in the Pleba\'nski complex. We will write down the operators most general as allowed by the relevant representation theory. The representation theory exercise is straightforward, and will not be spelled out. One calculates the different number of copies of one of the spaces in (\ref{complex-spaces}) in the tensor product of another space in (\ref{complex-spaces}) and the space $\Lambda^1$. One can then write down the relevant differential operators without difficulty. 

\subsection{Differential operators defined}

We now attempt to find the most general complex of the type
\be\label{complex-spaces}
\begin{tikzcd}
\Lambda^1 \arrow{r}{d_1} & S \arrow{r}{d_2} & E\times \Lambda^1 \arrow{r}{d_3} & E
\end{tikzcd}
\ee
The possible maps are dictated by the representation theory. They can be parametrised as
\be\label{d1-gen}
d_1 \xi = ( a_1 \partial^\mu \xi_\mu, a_2 \Sigma^{i\mu\nu} \partial_\mu \xi_\nu, 2 a_3 \partial_{\langle\mu} \xi_{\nu\rangle}) \in ( \Lambda^0, E, {\rm Sym}_0^2 T^*M) ,
\ee
where $a_1, a_2, a_3$ are arbitrary coefficients and
\be
\partial_{\langle\mu} \xi_{\nu\rangle}=
\frac{1}{2}\partial_{\mu} \xi_{\nu} + \frac{1}{2}\partial_{\nu} \xi_{\mu} - \frac{1}{4} g_{\mu\nu} \partial^\alpha \xi_\alpha
\ee
is the symmetric tracefree part. 
Here ${\rm Sym}_0^2 T^*M$ is the space of symmetric tracefree tensors. For the second map we have
\be\label{d2-gen}
d_2 \sigma = b_1 \Sigma^i_\mu{}^\alpha \partial_\alpha h + b_2 \partial_\mu h^i + b_3 \epsilon^{ijk} \Sigma^j_\mu{}^\alpha \partial_\alpha h^k + b_4 \Sigma^i_\mu{}^\rho \partial^\sigma \tilde{h}_{\rho\sigma} + b_5 \Sigma^{i\rho\sigma}  \partial_\rho \tilde{h}_{\mu\sigma} \in E\times\Lambda^1.
\ee
Again, $b_{1,2,3,4,5}$ are arbitrary coefficients. The last map is
\be\label{d3-gen}
d_3 a = c_1 \partial^\mu a_\mu^i + c_2 \epsilon^{ijk} \Sigma^{j\mu\nu} \partial_\mu a_\nu^k \in E.
\ee

\subsection{Compositions}

We want to the above to be a differential complex, so we want the compositions of the arrows to give zero. We have
\be
d_2 d_1 \xi = b_1 a_1 \Sigma^i_\mu{}^\rho \partial_\rho \partial^\sigma \xi_\sigma + b_2 a_2 \Sigma^{i\rho\sigma} \partial_\mu \partial_\rho \xi_\sigma + b_3 a_2 \epsilon^{ijk} \Sigma^j_\mu{}^\nu \partial_\nu \Sigma^{k\rho\sigma} \partial_\rho \xi_\sigma \\ \nonumber
+ b_4 a_3 \Sigma^i_\mu{}^\rho \partial^\sigma (\partial_\rho \xi_\sigma + \partial_\sigma \xi_\rho- \frac{1}{2} g_{\rho\sigma}\partial^\alpha \xi_\alpha ) 
+ b_5 a_3 \Sigma^{i\rho\sigma} \partial_\rho (\partial_\mu \xi_\sigma + \partial_\sigma \xi_\mu- \frac{1}{2} g_{\mu\sigma}\partial^\alpha \xi_\alpha).
\ee
We now use
\be
\epsilon^{ijk} \Sigma^j_{\mu\nu}  \Sigma^{k}_{\rho\sigma} = - g_{\mu\rho} \Sigma^i_{\nu\sigma} + g_{\nu\rho} \Sigma^i_{\mu\sigma} + g_{\mu\sigma} \Sigma^i_{\nu\rho} - g_{\nu\sigma} \Sigma^i_{\mu\rho}
\ee
to get
\be
d_2 d_1 \xi = (b_1 a_1+b_4 a_3-b_3 a_2 + \frac{b_5 a_3 - b_4 a_3}{2}) \Sigma^i_\mu{}^\rho \partial_\rho \partial^\sigma \xi_\sigma + (b_2 a_2+b_5 a_3-b_3 a_2) \Sigma^{i\rho\sigma} \partial_\mu \partial_\rho \xi_\sigma \\ \nonumber
+ (b_4 a_3+b_3 a_2) \Sigma^i_\mu{}^\rho \partial^\sigma  \partial_\sigma \xi_\rho.
\ee
This gives
\be\label{eqs-b1}
b_1 a_1+b_4 a_3-b_3 a_2+ \frac{b_5 a_3 - b_4 a_3}{2}=0, \qquad b_2 a_2+b_5 a_3-b_3 a_2=0, \qquad b_4 a_3+b_3 a_2=0.
\ee

The other composition is as follows
\be
d_3 d_2 \sigma = c_1 b_2 \partial^\mu \partial_\mu h^i + c_1 b_4 \Sigma^{i\mu\rho} \partial_\mu \partial^\sigma \tilde{h}_{\rho\sigma} + c_1 b_5 \Sigma^{i\rho\sigma} \partial_\rho \partial^\mu \tilde{h}_{\mu\sigma} \\ \nonumber
+ c_2 \epsilon^{ijk} \Sigma^{j\mu\nu} \partial_\mu ( b_1 \Sigma^k_\nu{}^\alpha \partial_\alpha h + b_2 \partial_\nu h^k + b_3 \epsilon^{klm} \Sigma^l_\nu{}^\alpha \partial_\alpha h^m 
+ b_4 \Sigma^k_\nu{}^\rho \partial^\sigma \tilde{h}_{\rho\sigma} + b_5 \Sigma^{k\rho\sigma} \partial_\rho \tilde{h}_{\nu\sigma}).
\ee
The first two terms in the second line do not contribute. The third term gives
\be
c_2 b_3 \epsilon^{ijk} \Sigma^{j\mu\nu} \partial_\mu \epsilon^{klm} \Sigma^l_\nu{}^\alpha \partial_\alpha h^m = c_2 b_3 (\delta^{il} \delta^{jm} - \delta^{im}\delta^{jl}) ( - \delta^{jl} g^{\mu\alpha} + \epsilon^{jls}\Sigma^{s\mu\alpha}) \partial_\mu\partial_\alpha h^m  \\ \nonumber
= 2 c_2 b_3 \partial^\mu \partial_\mu h^i.
\ee
The fourth term gives
\be
c_2 b_4 \epsilon^{ijk} \Sigma^{j\mu\nu}\Sigma^k_\nu{}^\rho  \partial_\mu \partial^\sigma \tilde{h}_{\rho\sigma}= 2 c_2 b_4 \Sigma^{i\mu\rho} \partial_\mu \partial^\sigma \tilde{h}_{\rho\sigma}.
\ee
The last term in the second line gives
\be\nonumber
c_2 b_5  \epsilon^{ijk} \Sigma^{j\mu\nu}  \Sigma^{k\rho\sigma} \partial_\mu \partial_\rho \tilde{h}_{\nu\sigma}=c_2 b_5 (- g^{\mu\rho} \Sigma^{i\nu\sigma} + g^{\nu\rho} \Sigma^{i\mu\sigma} + g^{\mu\sigma} \Sigma^{i\nu\rho} - g^{\nu\sigma} \Sigma^{i\mu\rho}) \partial_\mu \partial_\rho \tilde{h}_{\nu\sigma}= 0.
\ee
So, overall
\be
d_3 d_2 \sigma = (c_1 b_2 +2 c_2 b_3) \partial^\mu \partial_\mu h^i + (c_1 b_4 + c_1 b_5 + 2 c_2 b_4)\Sigma^{i\mu\rho} \partial_\mu \partial^\sigma \tilde{h}_{\rho\sigma}.
\ee
Requiring this to vanish gives
\be\label{eqs-b2}
c_1 b_2 +2 c_2 b_3=0, \qquad c_1 b_4 + c_1 b_5 + 2 c_2 b_4=0.
\ee
The system of equations (\ref{eqs-b1}), (\ref{eqs-b2}) can be interpreted as that for the constants $b_{1,2,3,4,5}$ in terms of $a_{1,2,3}, c_{1,2}$. The equations are not all independent, and not all 5 $b$'s can be solved for. Taking one of them as a parameter we have
\be
 b_2 = - \frac{2 a_1 b_1 c_2 }{a_2(c_1-c_2)}, \quad b_3 = \frac{ a_1 b_1 c_1 }{a_2(c_1-c_2)}, \quad b_4 = - \frac{a_1 b_1 c_1}{a_3(c_1-c_2)}, \quad b_5 = \frac{a_1 b_1 (c_1+2c_2)}{a_3 (c_1-c_2)}.
\ee
This means we can make the operator $d_2$ more explicit
\be\nonumber
d_2 \sigma = b (  \Sigma^i_\mu{}^\alpha \partial_\alpha h    - \frac{2 a_1 c_2 }{a_2(c_1-c_2)} \partial_\mu h^i +  \frac{ a_1 c_1 }{a_2(c_1-c_2)} \epsilon^{ijk} \Sigma^i_\mu{}^\alpha \partial_\alpha h^k \\ \nonumber
- \frac{a_1 c_1}{a_3(c_1-c_2)} \Sigma^i_\mu{}^\rho \partial^\sigma \tilde{h}_{\rho\sigma} + \frac{a_1 (c_1+2c_2)}{a_3 (c_1-c_2)} \Sigma^{i\rho\sigma}  \partial_\rho \tilde{h}_{\mu\sigma}),
\ee
where $b$ is an arbitrary parameter.

\subsection{Pleba\'nski case}

The Pleba\'nski case coefficients are
\be\label{Pleb-coeffs}
a_1=1, \quad a_2 = \frac{1}{4}, \quad a_3=\frac{1}{2}, \quad c_1=0, \quad c_2=1, \\ \nonumber
b_1=\frac{1}{4}, \quad b_2=2, \quad b_3=0, \quad b_4=0, \quad b_5=-1.
\ee
It can be checked that the compositions vanish with these coefficients, as they should. 

\subsection{Matching adjoints with inner products}

We now parametrise the adjoint operators as in (\ref{d1-adj-gen}), (\ref{d2-adj-gen}), (\ref{d3-adj-gen}). We will keep the inner product to be general, as in Section 5.2.
Let us relate the constants appearing in the adjoint operators to those in the inner products. We start with $d_1^*$. We have
\be
\langle \sigma, d_1 \xi\rangle = \langle d_1^* \sigma, \xi\rangle,
\ee
which gives
\be
\int_M \beta_1 a_1 h \partial^\mu \xi_\mu + \beta_2 a_2 h^i \Sigma^{i\mu\nu} \partial_\mu \xi_\nu + 2\beta_3 a_3 \tilde{h}^{\mu\nu} \partial_\mu \xi_\nu =\\ \nonumber
\int_M  \xi^\mu (a_1' \partial_\mu h + a_2' \Sigma^i_\mu{}^\nu \partial_\nu h^i + a_3' \partial^\nu \tilde{h}_{\mu\nu}).
\ee
So, we have
\be\label{a-prime}
\beta_1 a_1 = -  a_1', \qquad \beta_2 a_2 =  a_2', \qquad 2\beta_3 a_3 = -  a_3'.
\ee

Let us proceed similarly with $d_2^*$. We have
\be
\langle a, d_2 \sigma \rangle = \langle d_2^* a, \sigma \rangle,
\ee
which gives
\be\nonumber
\int_M \gamma_1 a^{i\mu} ( b_1 \Sigma^i_\mu{}^\alpha \partial_\alpha h + b_2 \partial_\mu h^i + b_3 \epsilon^{ijk} \Sigma^j_\mu{}^\alpha \partial_\alpha h^k + b_4 \Sigma^i_\mu{}^\rho \partial^\sigma \tilde{h}_{\rho\sigma} + b_5 \Sigma^{i\rho\sigma}  \partial_\rho \tilde{h}_{\mu\sigma}) \\ \nonumber
+ \gamma_2 \epsilon^{ijk} \Sigma^{i\mu\nu} a^j_\mu ( b_1 \Sigma^k_\nu{}^\alpha \partial_\alpha h + b_2 \partial_\nu h^k+ b_3 \epsilon^{klm} \Sigma^l_\nu{}^\alpha \partial_\alpha h^m + b_4 \Sigma^k_\nu{}^\rho \partial^\sigma \tilde{h}_{\rho\sigma} + b_5 \Sigma^{k\rho\sigma}  \partial_\rho \tilde{h}_{\nu\sigma}) = \\ \nonumber
\int_M \beta_1 b_1' h \Sigma^{i\mu\nu} \partial_\mu a_\nu^i + \beta_2 h^i (b_2' \partial^\mu a_\mu^i + b_3' \epsilon^{ijk} \Sigma^{j\mu\nu} \partial_\mu a_\nu^k) + \beta_3 \tilde{h}^{\mu\nu} (2b_4' \Sigma^i_{(\mu}{}^\alpha \partial_{\nu)}a_\alpha^i +  2b_5' \Sigma^i_{(\mu}{}^\alpha \partial_\alpha a_{\nu)}^i).
\ee

We first match the $h^i$ containing terms. Using the algebra of $\Sigma$'s gives
\be
\int_M ( \gamma_1 b_2 - 2 \gamma_2 b_3)a^{i\mu} \partial_\mu h^i + (\gamma_1 b_3 - \gamma_2 b_2 - \gamma_2 b_3) \epsilon^{ijk} a_\mu^i \Sigma^j_\mu{}^\alpha \partial_\alpha h^k  = \\ \nonumber
\int_M \beta_2 h^i (b_2' \partial^\mu a_\mu^i + b_3' \epsilon^{ijk} \Sigma^{j\mu\nu} \partial_\mu a_\nu^k) ,
\ee
which means
\be\label{b-prime-1}
\gamma_1 b_2 - 2 \gamma_2 b_3 = - \beta_2 b_2', \qquad \gamma_1 b_3 - \gamma_2 b_2 - \gamma_2 b_3 = - \beta_2 b_3'.
\ee

Let us now work out the $\tilde{h}_{\mu\nu}$ terms. The fourth term in the second line is
\be
\epsilon^{ijk} \Sigma^{i\mu\nu} a^j_\mu \Sigma^k_\nu{}^\rho \partial^\sigma \tilde{h}_{\rho\sigma}=\epsilon^{ijk}  \epsilon^{iks} \Sigma^{s\mu\rho} a^j_\mu\partial^\sigma \tilde{h}_{\rho\sigma}= - 2 \Sigma^{i\mu\rho} a^i_\mu\partial^\sigma \tilde{h}_{\rho\sigma}.
\ee
The last term in the second line
\be
\epsilon^{ijk} \Sigma^{i\mu\nu} a^j_\mu\Sigma^{k\rho\sigma}  \partial_\rho \tilde{h}_{\nu\sigma}= ( g^{\mu\rho} \Sigma^{i\nu\sigma} - g^{\nu\rho} \Sigma^{i\mu\sigma} - g^{\mu\sigma} \Sigma^{i\nu\rho} + g^{\nu\sigma} \Sigma^{i\mu\rho}) a^i_\mu \partial_\rho \tilde{h}_{\nu\sigma}= \\ \nonumber
-\Sigma^{i\mu\sigma}  a^i_\mu \partial^\nu \tilde{h}_{\nu\sigma} - \Sigma^{i\nu\rho} a^i_\mu \partial_\rho \tilde{h}_{\nu\mu} ,
\ee
where we have used that $\tilde{h}_{\mu\nu}$ is tracefree. This means that we have
\be
\int_M (\gamma_1 b_4 - 2  \gamma_2 b_4-\gamma_2 b_5 ) a^{i\mu} \Sigma^i_\mu{}^\rho \partial^\sigma \tilde{h}_{\rho\sigma} + (\gamma_1 b_5 + \gamma_2 b_5) \Sigma^{i\rho\sigma} a^{i\mu} \partial_\rho \tilde{h}_{\mu\sigma}  = \\ \nonumber
\int_M \beta_3 \tilde{h}^{\mu\nu} (2b_4' \Sigma^i_{(\mu}{}^\alpha \partial_{\nu)}a_\alpha^i +  2b_5' \Sigma^i_{(\mu}{}^\alpha \partial_\alpha a_{\nu)}^i).
\ee
This means that
\be\label{b-prime-2}
\gamma_1 b_4 - 2  \gamma_2 b_4-\gamma_2 b_5  = 2\beta_3 b_4', \qquad \gamma_1 b_5 + \gamma_2 b_5 = 2 \beta_3 b_5'.
\ee

We now match the $h$ terms
\be\nonumber
\int_M \gamma_1 b_1 a^{i\mu}  \Sigma^i_\mu{}^\alpha \partial_\alpha h    
+ \gamma_2 b_1 \epsilon^{ijk} \Sigma^{i\mu\nu} a^j_\mu \Sigma^k_\nu{}^\alpha \partial_\alpha h  = 
\int_M \beta_1 b_1' h \Sigma^{i\mu\nu} \partial_\mu a_\nu^i ,
\ee
which gives 
\be\label{b-prime-3}
b_1 (\gamma_1  - 2\gamma_2 ) = \beta_1 b_1' .
\ee

The calculation for the final adjoint is
\be
\langle \chi, d_3 a\rangle = \langle d_3^* \chi, a\rangle,
\ee
which reads
\be
\int_M  \chi^i (c_1 \partial^\mu a_\mu^i + c_2 \epsilon^{ijk} \Sigma^{j\mu\nu} \partial_\mu a_\nu^k ) =\\ \nonumber \int_M 
\gamma_1 a^{i\mu} (c_1' \partial_\mu \chi^i + c_2' \epsilon^{ijk} \Sigma^j_\mu{}^\alpha \partial_\alpha \chi^k) + \gamma_2 \epsilon^{ijk} \Sigma^{i\mu\nu} a_\mu^j (c_1' \partial_\nu \chi^k + c_2' \epsilon^{klm} \Sigma^l_\nu{}^\alpha \partial_\alpha \chi^m).
\ee
The last term in the second line simplifies
\be
\epsilon^{ijk} \Sigma^{i\mu\nu} a_\mu^j \epsilon^{klm} \Sigma^l_\nu{}^\alpha \partial_\alpha \chi^m = ( \delta^{il} \delta^{jm} - \delta^{im}\delta^{jl}) (-\delta^{il} g^{\mu\alpha} + \epsilon^{ils}\Sigma^{s\mu\alpha}) a_\mu^j \partial_\alpha \chi^m = \\ \nonumber
- 2 a^{i\mu} \partial_\mu \chi^i - \epsilon^{ijk} a_\mu^i \Sigma^{j\mu\alpha}  \partial_\alpha\chi^k.
\ee
This means that we have
\be
\int_M \chi^i (c_1 \partial^\mu a_\mu^i + c_2 \epsilon^{ijk} \Sigma^{j\mu\nu} \partial_\mu a_\nu^k ) =\\ \nonumber \int_M 
(\gamma_1 c_1' - 2 \gamma_2 c_2' ) a^{i\mu}  \partial_\mu \chi^i + (\gamma_1 c_2'  - \gamma_2 c_1' - \gamma_2 c_2' ) \epsilon^{ijk} a^{i\mu} \Sigma^j_\mu{}^\alpha \partial_\alpha \chi^k ,
\ee
and so
\be\label{c-prime}
- c_1 = \gamma_1 c_1' - 2 \gamma_2 c_2' , \qquad - c_2 = \gamma_1 c_2'  - \gamma_2 c_1' - \gamma_2 c_2' .
\ee
The equations we obtained give the coefficients appearing in the adjoints in terms of those in the operators and those in the inner products. 

\subsection{Composition of the adjoints}

Another useful calculation is that of the composition of the adjoints. We will compute the equations that guarantee that the compositions vanish. 
The object $d_2^* d_3^* \chi$ has three irreducible parts. The $h$ part equals
\be
b_1' \Sigma^{i\mu\nu} \partial_\mu ( c_1' \partial_\nu \chi^i + c_2' \epsilon^{ijk} \Sigma^j_\nu{}^\alpha \partial_\alpha\chi^k) =0.
\ee
The $h^i$ part is
\be
b_2' \partial^\mu (c_1' \partial_\mu \chi^i + c_2' \epsilon^{ijk} \Sigma^j_\mu{}^\alpha \partial_\alpha \chi^k) + b_3' \epsilon^{ijk} \Sigma^{j\mu\nu} \partial_\mu (c_1' \partial_\nu \chi^k + c_2' \epsilon^{klm} \Sigma^l_\nu{}^\alpha \partial_\alpha \chi^m) = \\ \nonumber
(b_2' c_1' + 2b_3' c_2' ) \partial^\mu  \partial_\mu \chi^i .
\ee
We thus want
\be
b_2' c_1' + 2b_3' c_2' =0.
\ee
The $\tilde{h}_{\mu\nu}$ part is the tracefree part of 
\be
2b_4' \Sigma^i_{(\mu}{}^\alpha \partial_{\nu)}( c_1' \partial_\alpha \chi^i + c_2' \epsilon^{ijk} \Sigma^j_\alpha{}^\beta \partial_\beta \chi^k) 
 +  2b_5' \Sigma^i_{(\mu}{}^\alpha \partial_\alpha ( c_1' \partial_{\nu)} \chi^i + c_2' \epsilon^{ijk} \Sigma^j_{\nu)}{}^\beta \partial_\beta \chi^k) .
 \ee
 The second term is
 \be
 \Sigma^i_{(\mu}{}^\alpha \partial_{\nu)}\epsilon^{ijk} \Sigma^j_\alpha{}^\beta \partial_\beta \chi^k = 2 \Sigma^i_{(\mu}{}^\alpha \partial_{\nu)} \partial_\alpha\chi^i.
 \ee
 The last term is
 \be
 \Sigma^i_{(\mu}{}^\alpha \partial_\alpha \epsilon^{ijk} \Sigma^j_{\nu)}{}^\beta \partial_\beta \chi^k = ( \delta^\alpha_{(\nu} \Sigma^i_{\mu)}{}^\beta + \delta^\beta_{(\mu} \Sigma^{i\alpha}{}_{\nu)}) \partial_\alpha\partial_\beta \chi^k=0.
 \ee
 So, overall we have
 \be
 (2b_4' c_1' + 4b_4' c_2' +2b_5' c_1') \Sigma^i_{(\mu}{}^\alpha \partial_{\nu)} \partial_\alpha \chi^i ,
 \ee
 and thus we want
 \be
 b_4' c_1' + 2b_4' c_2' +b_5' c_1'=0.
 \ee
 
 We now compute the other composition. We have
 \be
 d_1^* d_2^* a = a_1' b_1' \partial_\mu \Sigma^{i\rho\sigma} \partial_\rho a_\sigma^i + a_2' \Sigma^i_\mu{}^\nu \partial_\nu (b_2' \partial^\alpha a_\alpha^i + b_3' \epsilon^{ijk} \Sigma^{j\rho\sigma} \partial_\rho a_\sigma^k) \\ \nonumber
 + a_3' \partial^\nu ( b_4' \Sigma^i_{\mu}{}^\alpha \partial_{\nu}a_\alpha^i + b_4' \Sigma^i_{\nu}{}^\alpha \partial_{\mu}a_\alpha^i +  b_5' \Sigma^i_{\mu}{}^\alpha \partial_\alpha a_{\nu}^i+  b_5' \Sigma^i_{\nu}{}^\alpha \partial_\alpha a_{\mu}^i - \frac{b_4'-b_5'}{2} g_{\mu\nu} \Sigma^{i\rho\sigma}\partial_\rho a^i_\sigma).
 \ee
 The last term in the first line gives
 \be
 \Sigma^{i\mu\nu} \epsilon^{ijk} \Sigma^{j\rho\sigma} \partial_\nu \partial_\rho a_\sigma^k = ( - g^{\mu\rho} \Sigma^{i\nu\sigma} + g^{\nu\rho} \Sigma^{i\mu\sigma} + g^{\mu\sigma} \Sigma^{i\nu\rho} - g^{\nu\sigma} \Sigma^{i\mu\rho})  \partial_\nu \partial_\rho a_\sigma^i =\\ \nonumber 
 \Sigma^{i\mu\sigma} \partial^\alpha \partial_\alpha   a_\sigma^i
 -  \Sigma^{i\rho\sigma}  \partial^\mu \partial_\rho a_\sigma^i -   \Sigma^{i\mu\rho}  \partial_\rho  \partial^\sigma a_\sigma^i.
 \ee
 This gives
 \be\nonumber
 d_1^* d_2^* a = (a_1' b_1' + a_3' b_4' -  a_2' b_3'  - \frac{a_3'(b_4'-b_5')}{2}) \partial_\mu \Sigma^{i\rho\sigma} \partial_\rho a_\sigma^i + (a_2' b_2'  +  a_3' b_5' -   a_2' b_3' ) \Sigma^i_\mu{}^\nu \partial_\nu  \partial^\alpha a_\alpha^i \\ \nonumber
+ (a_2' b_3' + a_3' b_4' ) \Sigma^{i}_\mu{}^{\sigma} \partial^\alpha \partial_\alpha   a_\sigma^i .
 \ee
We thus want
\be
a_1' b_1' + a_3' b_4' -  a_2' b_3'  - \frac{a_3'(b_4'-b_5')}{2} =0, \qquad a_2' b_2'  +  a_3' b_5' -   a_2' b_3' =0, \qquad a_2' b_3' + a_3' b_4' =0.
\ee
We have checked that all these equations hold provided the primed coefficients are as for the adjoint operators, and provided that squaring relations (\ref{eqs-b1}), (\ref{eqs-b2})  hold.

\subsection{Computation of the square}
\label{sec:square}

We now form an elliptic operator
\be\label{D}
(S,E)\ni (\sigma,\chi) \to D(\sigma,\chi) = (d_1^*\sigma, d_2 \sigma + d_3^*\chi)\in (\Lambda^1,E\times\Lambda^1).
\ee
Its adjoint is given by 
\be\label{D*}
(\Lambda^1,E\times\Lambda^1) \ni (\xi,a) \to D^*(\xi,a) = ( d_1\xi + d_2^* a, d_3 a) \in (S,E).
\ee
The composition of these two operators is
\be
D^* D (\sigma,\chi) = ( d_1 d_1^*\sigma + d_2^* (d_2 \sigma + d_3^*\chi), d_3 (d_2 \sigma + d_3^*\chi))=( (d_1 d_1^* +d_2^* d_2 ) \sigma, d_3  d_3^*\chi).
\ee
We want to compute all the operators appearing here. 

We start with $d_1 d_1^*$. We have the following $h$ component of $d_1 d_1^*\sigma$
\be
 a_1 \partial^\mu (a_1' \partial_\mu h + a_2' \Sigma^i_\mu{}^\nu \partial_\nu h^i + a_3' \partial^\nu \tilde{h}_{\mu\nu})= a_1 a_1' \partial^\mu \partial_\mu h + a_1 a_3' \partial^\mu \partial^\nu \tilde{h}_{\mu\nu}.
 \ee
 For the $h^i$ component we have
 \be
 a_2 \Sigma^{i\mu\nu} \partial_\mu (a_1' \partial_\nu h + a_2' \Sigma^j_\nu{}^\alpha \partial_\alpha h^j + a_3' \partial^\alpha \tilde{h}_{\nu\alpha})=
 - a_2 a_2' \partial^\mu \partial_\mu h^i+ a_2 a_3' \Sigma^{i\mu\nu} \partial_\mu \partial^\alpha \tilde{h}_{\nu\alpha}.
 \ee
 For the $\tilde{h}_{\mu\nu}$ component we have
 \be
 2 a_3 \partial_{(\mu} (a_1' \partial_{\nu)} h + a_2' \Sigma^i_{\nu)}{}^\alpha \partial_\alpha h^i + a_3' \partial^\alpha \tilde{h}_{\nu)\alpha})=\\ \nonumber
 2a_3 a_1'  \partial_{\mu} \partial_{\nu} h + 2a_3 a_2'  \Sigma^i_{(\mu}{}^\alpha \partial_{\nu)} \partial_\alpha h^i + 2 a_3 a_3' \partial_{(\mu}\partial^\alpha \tilde{h}_{\nu)\alpha},
 \ee
 where it is understood that the tracefree part of the result is taken. 
 We now compute $d_2^* d_2  \sigma$. For the $h$ component we have
 \be
  b_1' \Sigma^{i\mu\nu} \partial_\mu (b_1 \Sigma^i_\nu{}^\alpha \partial_\alpha h + b_2 \partial_\nu h^i + b_3 \epsilon^{ijk} \Sigma^j_\nu{}^\alpha \partial_\alpha h^k + b_4 \Sigma^i_\nu{}^\rho \partial^\sigma \tilde{h}_{\rho\sigma} + b_5 \Sigma^{i\rho\sigma}  \partial_\rho \tilde{h}_{\nu\sigma})= \\ \nonumber
  - 3 b_1 b_1' \partial^\mu \partial_\mu h - 3 b_4 b_1' \partial^\mu \partial^\nu \tilde{h}_{\mu\nu} + b_5 b_1' ( g^{\mu\rho} g^{\nu\sigma}- g^{\mu\sigma} g^{\nu\rho} + \epsilon^{\mu\nu\rho\sigma})
  \partial_\mu \partial_\rho \tilde{h}_{\nu\sigma} = \\ \nonumber
  ( - 3 b_1 b_1' + b_5 b_1' ) \partial^\mu \partial_\mu h - (3 b_4 b_1'  + b_5 b_1' )\partial^\mu \partial^\nu \tilde{h}_{\mu\nu}.
  \ee
  For the $h^i$ component we have
  \be\nonumber
  b_2' \partial^\mu (b_1 \Sigma^i_\mu{}^\alpha \partial_\alpha h + b_2 \partial_\mu h^i + b_3 \epsilon^{ijk} \Sigma^j_\mu{}^\alpha \partial_\alpha h^k + b_4 \Sigma^i_\mu{}^\rho \partial^\sigma \tilde{h}_{\rho\sigma} + b_5 \Sigma^{i\rho\sigma}  \partial_\rho \tilde{h}_{\mu\sigma}) \\ \nonumber
  + b_3' \epsilon^{ijk} \Sigma^{j\mu\nu} \partial_\mu (b_1 \Sigma^k_\nu{}^\alpha \partial_\alpha h + b_2 \partial_\nu h^k + b_3 \epsilon^{klm} \Sigma^l_\nu{}^\alpha \partial_\alpha h^m + b_4 \Sigma^k_\nu{}^\rho \partial^\sigma \tilde{h}_{\rho\sigma} + b_5 \Sigma^{k\rho\sigma}  \partial_\rho \tilde{h}_{\nu\sigma}) = \\ \nonumber
  b_2  b_2' \partial^\mu \partial_\mu h^i+  (b_4 b_2' +  b_5 b_2' )  \Sigma^i_\mu{}^\rho \partial^\mu \partial^\sigma \tilde{h}_{\rho\sigma} 
  + b_3 b_3' (\delta^{il}\delta^{jm} - \delta^{im}\delta^{jl}) ( - \delta^{jl} g^{\mu\alpha} + \epsilon^{jls}\Sigma^{s\mu\alpha}) 
   \partial_\mu \partial_\alpha h^m \\ \nonumber
   + 2b_4 b_3' \Sigma^{i\mu\rho} \partial_\mu \partial^\sigma \tilde{h}_{\rho\sigma}
     + b_5 b_3' (- g^{\mu\rho} \Sigma^{i\nu\sigma} + g^{\nu\rho} \Sigma^{i\mu\sigma} + g^{\mu\sigma} \Sigma^{i\nu\rho} - g^{\nu\sigma} \Sigma^{i\mu\rho})   \partial_\mu \partial_\rho \tilde{h}_{\nu\sigma} = \\ \nonumber
      (b_2  b_2' + 2 b_3 b_3' )\partial^\mu \partial_\mu h^i+  (b_4 b_2' +  b_5 b_2' + 2b_4 b_3' )  \Sigma^i_\mu{}^\rho \partial^\mu \partial^\sigma \tilde{h}_{\rho\sigma} .
    \ee
  For the $\tilde{h}_{\mu\nu}$ component we need the tracefree part of
  \be\nonumber
  2b_4' \Sigma^i_{(\mu}{}^\alpha \partial_{\nu)}(b_1 \Sigma^i_\alpha{}^\beta \partial_\beta h + b_2 \partial_\alpha h^i + b_3 \epsilon^{ijk} \Sigma^j_\alpha{}^\beta \partial_\beta h^k + b_4 \Sigma^i_\alpha{}^\rho \partial^\sigma \tilde{h}_{\rho\sigma} + b_5 \Sigma^{i\rho\sigma}  \partial_\rho \tilde{h}_{\alpha\sigma}) \\ \nonumber
  +  2b_5' \Sigma^i_{(\mu}{}^\alpha \partial_\alpha (b_1 \Sigma^i_{\nu)}{}^\beta \partial_\beta h + b_2 \partial_{\nu)} h^i + b_3 \epsilon^{ijk} \Sigma^j_{\nu)}{}^\beta \partial_\beta h^k + b_4 \Sigma^i_{\nu)}{}^\rho \partial^\sigma \tilde{h}_{\rho\sigma} + b_5 \Sigma^{i\rho\sigma}  \partial_\rho \tilde{h}_{\nu)\sigma}) = \\ \nonumber
 -6 b_1 b_4' \partial_{(\mu} \partial_{\nu)} h + (2b_2 b_4' + 4b_3 b_4' ) \Sigma^i_{(\mu}{}^\alpha \partial_{\nu)} \partial_\alpha h^i
  -(6 b_4 b_4'  + 2b_5 b_4' ) \partial_{(\mu} \partial^\sigma \tilde{h}_{\nu)\sigma} \\ \nonumber
  +2 b_1b_5' ( g_{\mu\nu} \partial^\alpha\partial_\alpha h - \partial_{(\mu} \partial_{\nu)} h) 
  +  2b_2 b_5' \Sigma^i_{(\mu}{}^\alpha \partial_\alpha \partial_{\nu)} h^i 
  +  2b_3 b_5'  ( - g_{\mu\nu} \Sigma^{i\alpha\beta} + g_{(\mu}^\alpha  \Sigma^{i}_{\nu)}{}^{\beta} - g_{(\mu}^\beta \Sigma^{i}_{\nu)}{}^\alpha) \partial_\alpha  \partial_\beta h^i
  \\ \nonumber
  +  2b_4 b_5' ( g_{\mu\nu} \partial^\rho\partial^\sigma \tilde{h}_{\rho\sigma} - \partial_{(\mu} \partial^\rho \tilde{h}_{\nu)\rho}) 
  +  2 b_5 b_5' (\partial_{(\mu} \partial^\rho \tilde{h}_{\nu)\rho} - \partial^\alpha \partial_\alpha \tilde{h}_{\mu\nu}) =  \\ \nonumber
  +( -6 b_1 b_4'  - 2 b_1b_5' ) \partial_{(\mu} \partial_{\nu)} h  +2 b_1b_5' g_{\mu\nu} \partial^\alpha\partial_\alpha h + (2b_2 b_4' + 4b_3 b_4' +  2b_2 b_5' ) \Sigma^i_{(\mu}{}^\alpha \partial_{\nu)} \partial_\alpha h^i \\ \nonumber
  +(-6 b_4 b_4'  - 2b_5 b_4' -  2b_4 b_5'+  2 b_5 b_5') \partial_{(\mu} \partial^\sigma \tilde{h}_{\nu)\sigma} 
  +  2b_4 b_5'  g_{\mu\nu} \partial^\rho\partial^\sigma \tilde{h}_{\rho\sigma} 
  -  2 b_5 b_5'  \partial^\alpha \partial_\alpha \tilde{h}_{\mu\nu}.
     \ee
Here we have used that $\tilde{h}_{\mu\nu}$ is tracefree. Dropping the trace parts, and collecting all the terms, we have the following expression for the $\tilde{h}_{\mu\nu}$ part of $(d_1 d_1^* +d_2^* d_2 ) \sigma$
 \be
 ( -6 b_1 b_4'  - 2 b_1b_5' + 2a_3 a_1' ) \partial_{\langle\mu} \partial_{\nu\rangle} h  + (2b_2 b_4' + 4b_3 b_4' +  2b_2 b_5' + 2a_3 a_2' ) \Sigma^i_{\langle\mu}{}^\alpha \partial_{\nu\rangle} \partial_\alpha h^i \\ \nonumber
  +(-6 b_4 b_4'  - 2b_5 b_4' -  2b_4 b_5'+  2 b_5 b_5'+ 2 a_3 a_3' ) \partial_{\langle\mu} \partial^\sigma \tilde{h}_{\nu\rangle\sigma} 
  -  2 b_5 b_5'  \partial^\alpha \partial_\alpha \tilde{h}_{\mu\nu}.
  \ee
  For the $h^i$ part we have
  \be
     (b_2  b_2' + 2 b_3 b_3'  - a_2 a_2' )\partial^\mu \partial_\mu h^i+  (b_4 b_2' +  b_5 b_2' + 2b_4 b_3' + a_2 a_3'  ) \Sigma^i_\mu{}^\rho \partial^\mu \partial^\sigma \tilde{h}_{\rho\sigma}.
    \ee
For the trace part $h$ we have 
    \be
     ( - 3 b_1 b_1' + b_5 b_1' + a_1 a_1' ) \partial^\mu \partial_\mu h + (-3 b_4 b_1'  - b_5 b_1' + a_1 a_3' )\partial^\mu \partial^\nu \tilde{h}_{\mu\nu}.  
      \ee
      
  We also need to compute the operator $d_3 d_3^* \chi$. We have
  \be
  d_3 d_3^* \chi = c_1 \partial^\mu (c_1' \partial_\mu \chi^i + c_2' \epsilon^{ijk} \Sigma^j_\mu{}^\alpha \partial_\alpha \chi^k) 
  + c_2 \epsilon^{ijk} \Sigma^{j\mu\nu} \partial_\mu (c_1' \partial_\nu \chi^k + c_2' \epsilon^{klm} \Sigma^l_\nu{}^\alpha \partial_\alpha \chi^m) = \\ \nonumber
 ( c_1 c_1' + 2 c_2 c_2')\partial^\mu \partial_\mu \chi^i.
  \ee
      
\subsection{Imposing the \texorpdfstring{$\Delta$}{} conditions}
      
 We want the $D^* D$ operator to be a multiple of the $\Delta=\partial^\alpha \partial_\alpha$ operator. This means we want to have the following quantities vanishing
 \be\label{D2-eqs}
 -6 b_1 b_4' - 2 b_1b_5' + 2a_3 a_1' =0, \quad 2b_2 b_4' + 4b_3 b_4' +  2b_2 b_5' + 2a_3 a_2' =0, \\ \nonumber
  -6 b_4 b_4'  - 2b_5 b_4' -  2b_4 b_5'+  2 b_5 b_5'+ 2 a_3 a_3' =0, 
 \\ \nonumber
 b_4 b_2' +  b_5 b_2' + 2b_4 b_3' + a_2 a_3' =0, \quad -3 b_4 b_1'  - b_5 b_1' + a_1 a_3'=0.
 \ee
 When these conditions are satisfied, the parts of $D^* D$ become
 \be\nonumber
 (d_1 d_1^* +d_2^* d_2 ) \sigma =  ( - 3 b_1 b_1' + b_5 b_1' + a_1 a_1' ) \partial^\mu \partial_\mu h +(b_2  b_2' + 2 b_3 b_3'  - a_2 a_2' )\partial^\mu \partial_\mu h^i -  2 b_5 b_5'  \partial^\alpha \partial_\alpha h_{\mu\nu}, \\ \nonumber
 d_3 d_3^* \chi= ( c_1 c_1' + 2 c_2 c_2')\partial^\mu \partial_\mu \chi^i.
 \ee
Solving the conditions (\ref{D2-eqs}) for the operators $d_1, d_2, d_3$ of the Pleba\'nski complex leads to the choice of the inner products as in (\ref{inner-prod-1}). As we have already described, the resulting in this case operator $d_2^*$ is not the one relevant for the linearised Einstein equations. This means we need a different operator $D$. We will build it in the next section, based on the general computations that were carried out in this section.  
 
\section{Twisted operator}

\subsection{Construction of a more general operator \texorpdfstring{$\tilde{D}$}{}}

We know that the operators $D,D^*$  (\ref{D}), (\ref{D*}) constructed from the Pleba\'nski complex do not satisfy $D^* D\sim\Delta$. We would now like to use the operators of the Pleba\'nski complex as building blocks of more general operators $\tilde{D},\tilde{D}^*$. Thus, we consider
\be\label{tilde-D}
\tilde{D}(\sigma,\chi) = (\tilde{d}_1^* \sigma + \tilde{d}_4 \chi, \tilde{d}_2 \sigma + \tilde{d}_3^*\chi),
\ee
Here all operators are general, of the type considered in the previous section, see (\ref{d1-gen}), (\ref{d2-gen}), (\ref{d3-gen}), and we have introduced a new operator $d_4:E\to \Lambda^1$
\be
\tilde{d}_4 \chi=f \Sigma^i_\mu{}^\alpha \partial_\alpha \chi^i .
 \ee
 Let us also introduce
\be
\tilde{D}^*(\xi, a) = (\tilde{d}_1 \xi + \tilde{d}_2^* a, \tilde{d}_4^* \xi + \tilde{d}_3 a).
\ee
We do not yet assume that the operators appearing here are adjoints of those in (\ref{tilde-D}) with respect to some inner product. For now these are just the most general operators of the type (\ref{d1-adj-gen}), (\ref{d2-adj-gen}), (\ref{d3-adj-gen}) and we introduced
\be
\tilde{d}_4^* \xi = f' \Sigma^{i\mu\nu} \partial_\mu\xi_\nu.
\ee

The composition $\tilde{D}^* \tilde{D}$ is 
\be
\tilde{D}^* \tilde{D}(\sigma,\chi) = (\tilde{d}_1 (\tilde{d}_1^* \sigma + \tilde{d}_4 \chi) + \tilde{d}_2^* (\tilde{d}_2 \sigma + \tilde{d}_3^*\chi), \tilde{d}_4^* (\tilde{d}_1^* \sigma + \tilde{d}_4 \chi)+ \tilde{d}_3 (\tilde{d}_2 \sigma + \tilde{d}_3^*\chi))
\ee
Most of the operators appearing here were already evaluated in section \ref{sec:square}. Let us compute the parts that have not been previously computed.

\subsection{\texorpdfstring{Computation of $\tilde{d}_4^* \tilde{d}_1^* +\tilde{d}_3 \tilde{d}_2$}{}}

We have
\be\nonumber
\tilde{d}_4^* \tilde{d}_1^* +\tilde{d}_3 \tilde{d}_2 =f' \Sigma^{i\mu\nu} \partial_\mu(a_1'  \partial_\nu h + a_2' \Sigma^j_\nu{}^\alpha \partial_\alpha h^j +a_3' \partial^\alpha \tilde{h}_{\nu\alpha}) \\ \nonumber
+ c_1 \partial^\mu ( b_1 \Sigma^i_\mu{}^\alpha \partial_\alpha h + b_2 \partial_\mu h^i + b_3 \epsilon^{ijk} \Sigma^j_\mu{}^\alpha \partial_\alpha h^k + b_4 \Sigma^i_\mu{}^\rho \partial^\sigma \tilde{h}_{\rho\sigma} + b_5 \Sigma^{i\rho\sigma}  \partial_\rho \tilde{h}_{\mu\sigma} ) \\ \nonumber
+  c_2 \epsilon^{ijk} \Sigma^{j\mu\nu} \partial_\mu ( b_1 \Sigma^k_\nu{}^\alpha \partial_\alpha h + b_2 \partial_\nu h^k + b_3 \epsilon^{klm} \Sigma^l_\nu{}^\alpha \partial_\alpha h^m + b_4 \Sigma^k_\nu{}^\rho \partial^\sigma \tilde{h}_{\rho\sigma} + b_5 \Sigma^{k\rho\sigma}  \partial_\rho \tilde{h}_{\nu\sigma} ) =  \\ \nonumber
- f' a_2' \partial^\alpha \partial_\alpha h^i
+ f' a_3'\Sigma^{i\mu\nu} \partial_\mu \partial^\alpha \tilde{h}_{\nu\alpha}  + c_1 b_2 \partial^\mu \partial_\mu h^i + c_1 (b_4 +b_5) \Sigma^{i\mu\rho} \partial_\mu \partial^\sigma \tilde{h}_{\rho\sigma} 
\\ \nonumber
+ c_2 b_3 \epsilon^{ijk} \Sigma^{j\mu\nu}  \epsilon^{klm} \Sigma^l_\nu{}^\alpha \partial_\mu \partial_\alpha h^m + c_2 b_4 \epsilon^{ijk} \Sigma^{j\mu\nu}
\Sigma^k_\nu{}^\rho \partial_\mu \partial^\sigma \tilde{h}_{\rho\sigma} + c_2 b_5 \epsilon^{ijk} \Sigma^{j\mu\nu} \Sigma^{k\rho\sigma}  \partial_\mu \partial_\rho \tilde{h}_{\nu\sigma}.
\ee
The terms in the last line compute as follows
\be
\epsilon^{ijk} \Sigma^{j\mu\nu}  \epsilon^{klm} \Sigma^l_\nu{}^\alpha \partial_\mu \partial_\alpha h^m= (\delta^{il} \delta^{jm} - \delta^{im}\delta^{jl}) ( - \delta^{jl} g^{\mu\alpha} + \epsilon^{jls}\Sigma^{s\mu\alpha})  \partial_\mu \partial_\alpha h^m= 2 \partial^\mu \partial_\mu h^i, \\ \nonumber
\epsilon^{ijk} \Sigma^{j\mu\nu}\Sigma^k_\nu{}^\rho \partial_\mu \partial^\sigma \tilde{h}_{\rho\sigma}= 2 \Sigma^{i\mu\rho} \partial_\mu \partial^\sigma \tilde{h}_{\rho\sigma}, \\ \nonumber
\epsilon^{ijk} \Sigma^{j\mu\nu} \Sigma^{k\rho\sigma}  \partial_\mu \partial_\rho \tilde{h}_{\nu\sigma} = (- g^{\mu\rho} \Sigma^{i\nu\sigma} + g^{\nu\rho} \Sigma^{i\mu\sigma} + g^{\mu\sigma} \Sigma^{i\nu\rho} - g^{\nu\sigma} \Sigma^{i\mu\rho})\partial_\mu \partial_\rho \tilde{h}_{\nu\sigma}= 0.
\ee
This gives
\be\nonumber
\tilde{d}_4^* \tilde{d}_1^* +\tilde{d}_3 \tilde{d}_2 =(- f' a_2' +2 c_2 b_3+ c_1 b_2) \partial^\alpha \partial_\alpha h^i
+ (f' a_3' + c_1 (b_4 +b_5) + 2 c_2 b_4)\Sigma^{i\mu\nu} \partial_\mu \partial^\alpha \tilde{h}_{\nu\alpha}.
\ee
In order for our operator $\tilde{D}$ to square to a multiple of the box operator the second term here must vanish and so we impose
\be\label{f-eqn}
f' a_3' + c_1 (b_4 +b_5) + 2 c_2 b_4=0.
\ee
This can be viewed as an equation giving $f'$ provided all other parameters are known.

\subsection{Finding \texorpdfstring{$\tilde{D},\tilde{D}^*$}{} that encode Einstein equations mod gauge}

We now look for the solution of the system of equations $\tilde{D}^* \tilde{D}\sim\Delta$ assuming that $\tilde{d}_2=d_2$ and $\tilde{d}_2^* = d_2^*$ of the Pleba\'nski complex. These are the two operators whose composition $d_2^* d_2$ encodes the linearised Einstein equations, and we would like to keep this part of $\tilde{D}^*\tilde{D}$ intact. All other entries in $\tilde{D}^*\tilde{D}$ are related to gauge, and we are at will to choose them as convenient. Choosing the coefficients $b,b'$ as they are in the Pleba\'nski complex case, and substituting these into (\ref{D2-eqs}), we get 
\be
a_2= a_2'=0, \quad a_1=-2a_3, \quad a_1'= \frac{1}{4a_3}, \quad a_3'=-\frac{1}{a_3}, \quad f' = - c_1 a_3.
\ee
This gives
\be
\tilde{d}_1 \xi = a_3 ( - 2\partial^\mu \xi_\mu, 0, 2\partial_{\langle\mu} \xi_{\nu\rangle} ), \\ \nonumber
\tilde{d}_1^* \sigma =  \frac{1}{4a_3} \partial_\mu h - \frac{1}{a_3} \partial^\nu \tilde{h}_{\mu\nu}.
\ee
We can choose $a_3$ in such a way that $\tilde{d}_1^*$ is the adjoint of $\tilde{d}_1$ with respect to the inner products (\ref{inner-prod-pleb}). We have
\be
(\tilde{d}_1\xi,\sigma) = \int -\frac{a_3}{2} (\partial^\mu \xi_\mu) h + 2a_3 \tilde{h}^{\mu\nu} \partial_\mu \xi_\nu = \int ( \frac{a_3}{2} \partial_\mu h - 2a_3 \partial^\nu \tilde{h}_{\mu\nu}) \xi^\mu.
\ee
So, the adjoint of $\tilde{d}_1$ matches $\tilde{d}_1^*$ if 
\be
\frac{a_3}{2}= \frac{1}{4a_3}, \qquad - 2a_3 = - \frac{1}{a_3} \qquad \Rightarrow \qquad a_3=\frac{1}{\sqrt{2}}.
\ee
Thus, finally
\be\label{tilde-d1-d1*}
\tilde{d}_1 \xi = \sqrt{2} ( - \partial^\mu \xi_\mu, 0, \partial_{\langle\mu} \xi_{\nu\rangle} ), \\ \nonumber
\tilde{d}_1^* \sigma =  \sqrt{2}( \frac{1}{4} \partial_\mu h -  \partial^\nu \tilde{h}_{\mu\nu}).
\ee

As a check we compute
\be
\tilde{d}_1 \tilde{d}_1^* \sigma = 2( - \frac{1}{4}\partial^\mu \partial_\mu h + \partial^\mu \partial^\nu \tilde{h}_{\mu\nu}, 0, \frac{1}{4} \partial_{\langle\mu} \partial_{\nu\rangle} h -  \partial_{\langle \mu} \partial^\alpha \tilde{h}_{\nu\rangle\alpha}).
\ee
Recalling (\ref{d2}) we have
\be
(d_2^* d_2 + \tilde{d}_1 \tilde{d}_1^*)\sigma = (\partial^\mu \partial_\mu h, 0, - \partial^\alpha \partial_\alpha \tilde{h}_{\mu\nu}),
\ee
which is as desired. Note that the composition $d_2 \tilde{d}_1$ does not vanish, so these operators do not form a complex. However, as we will now see, the tilded operators can be rewritten in terms of those of Pleba\'nski complex. 

\subsection{Rewriting} 

We would like to rewrite the operators $\tilde{d}_1, \tilde{d}_1^*$ in terms of those appearing in the Pleba\'nski complex. For this, we introduce two linear operators 
\be\label{Phi}
\Phi: E\otimes \Lambda^1 \to \Lambda^1, \qquad \Phi(a)_\mu = \Sigma^i_{\mu}{}^\alpha a^i_\alpha,
\ee
and
\be
\Phi^*: \Lambda^1 \to E\otimes \Lambda^1, \qquad \Phi^*(\xi) = \frac{1}{2 \gamma_2 - \gamma_1} \Sigma^i_\mu{}^\alpha \xi_\alpha
\ee
such that
\be
\langle \xi, \Phi(a) \rangle = \langle a^i, \Phi^*(\xi)^i \rangle
\ee
where $\Phi^*$ is the adjoint of $\Phi$.
Note that for the case of Pleba\'nski $\gamma_1 = 0, \gamma_2 = 1$ so
\be
\Phi^*(\xi)^i_\mu = \frac{1}{2} \Sigma^i_\mu{}^\alpha \xi_\alpha.
\ee
We then compute
\be\label{Phi-d2}
\Phi (d_2 \sigma)= \Sigma^i_\mu{}^\alpha( \frac{1}{4} \Sigma^i_\alpha{}^\nu \partial_\nu h +2 \partial_\alpha h^i -  \Sigma^{i\rho\sigma} \partial_\rho \tilde{h}_{\alpha\sigma}) =
 \\ \nonumber
-\frac{3}{4} \partial_\mu h + 2 \Sigma^i_\mu{}^\alpha\partial_\alpha h^i - ( \delta_\mu{}^\rho g^{\alpha\sigma} - \delta_\mu{}^\sigma g^{\alpha\rho}) \partial_\rho \tilde{h}_{\alpha\sigma} =
- \frac{3}{4} \partial_\mu h + \partial^\nu \tilde{h}_{\mu\nu} + 2 \Sigma^i_\mu{}^\alpha\partial_\alpha h^i .
\ee
It is now clear that there is a linear combination of $d_1^*\sigma$ and $\Sigma^i_\mu{}^\alpha (d_2 \sigma)_\alpha^i$ in which the $h^i$ term cancels
\be\label{tilde-d1-star}
(d_1^*\sigma)_\mu - \Phi (d_2 \sigma)_\mu = -\frac{1}{4}\partial_\mu h -  \partial^\nu \tilde{h}_{\mu\nu} + \frac{3}{4} \partial_\mu h - \partial^\nu \tilde{h}_{\mu\nu}
=  \frac{1}{2} \partial_\mu h-2 \partial^\nu \tilde{h}_{\mu\nu} =  \sqrt{2} \tilde{d}_1^* \sigma.
\ee
Notably, we were able to rewrite the operator  $\tilde{d}_1^*$ that we deduced is necessary for $\tilde{D}^* \tilde{D}\sim\Delta$ in terms of the original Pleba\'nski operators. In fact, the inner products (\ref{inner-prod-pleb}) were chosen as they are so that this becomes possible. 

We also have
\be
2 d_2^*( \Phi^*(\xi)) = (-2\Sigma^{i\mu\nu}  \Sigma^i_\nu{}^\alpha \partial_\mu \xi_\alpha, \frac{1}{4} \epsilon^{ijk} \Sigma^{j\mu\nu}  \Sigma^k_\nu{}^\alpha \partial_\mu \xi_\alpha, \Sigma^i_{\langle\mu|}{}^\alpha \Sigma^i_\alpha{}^\beta \partial_{|\nu\rangle} \xi_\beta-\Sigma^i_{\langle\mu|}{}^\alpha \Sigma^i_{|\nu\rangle}{}^\beta \partial_\alpha  \xi_\beta ) \\ \nonumber = 
 (  6 \partial^\mu \xi_\mu, \frac{1}{2} \Sigma^{i\mu\nu} \partial_\mu  \xi_\nu, -2 \partial_{\langle\mu}\xi_{\nu\rangle} ).
\ee
This means we have
\be
d_1\xi - d_2^*( \Phi^*(\xi)) = (-2 \partial^\mu \xi_\mu , 0, 2 \partial_{\langle\mu}\xi_{\nu\rangle} ) = \sqrt{2} \tilde{d}_1 \xi.
\ee
Again, we were able to rewrite the tilded operator in terms of the original Pleba\'nski complex untilded ones. This means that the operators $\tilde{D},\tilde{D}^*$ start to take a more concrete form
\be
\tilde{D}(\sigma,\chi) = (\frac{1}{\sqrt{2}}(d_1^*\sigma - \Phi(d_2\sigma)) + \tilde{d}_4 \chi, d_2\sigma + \tilde{d}_3^* \chi), \\ \nonumber
\tilde{D}^*(\xi,a) = ( \frac{1}{\sqrt{2}} (d_1 \xi -  d_2^* (\Phi^*(\xi))) + d_2^* a, \tilde{d}_3 a + \tilde{d}_4^* \xi).
\ee
It remains to chose the operators $\tilde{d}_3, \tilde{d}_4$ and $\tilde{d}_3^*, \tilde{d}_4^*$ in the most convenient way.

\subsection{Requiring the adjoints}

We have imposed the condition that $\tilde{d}_1, \tilde{d}_1^*$ are adjoints of each other. Let us demand that the same holds true for the operators $\tilde{d}_3, \tilde{d}_4$ and $\tilde{d}_3^*, \tilde{d}_4^*$. The general operator $\tilde{d}_3$ is parametrised
\be
\tilde{d}_3 a = c_1 \partial^\mu a_\mu^i + c_2 \epsilon^{ijk} \Sigma^{j\mu\nu} \partial_\mu a_\nu^k.
\ee
In the inner products (\ref{inner-prod-pleb}) the coefficients of its adjoint are $c_1'=c_2 - c_1/2, c_2'   = c_1/2$, so that 
\be
\tilde{d}_3^* \chi = ( c_2 - \frac{c_1}{2}) \partial_\mu \chi^i + \frac{c_1}{2} \epsilon^{ijk} \Sigma^j_\mu{}^\alpha \partial_\alpha \chi^k.
\ee
For the operators $\tilde{d}_4, \tilde{d}_4^*$ we have
\be
(\xi, \tilde{d}_4 \chi) = \int f \xi^\mu \Sigma^i_\mu{}^\alpha \partial_\alpha \chi^i = (\tilde{d}_4^* \xi, \chi),
\ee
where 
\be
\tilde{d}_4^* \xi = f \Sigma^{i\mu\nu} \partial_\mu \xi_\nu,
\ee
so that $f'=f$. 

\subsection{Additional computations}

We have
\be\nonumber
(\tilde{d}_4^* \tilde{d}_1^* +\tilde{d}_3 d_2)\sigma =f' \Sigma^{i\mu\nu} \partial_\mu \sqrt{2} (\frac{1}{4} \partial_\nu h -  \partial^\alpha \tilde{h}_{\nu\alpha}) 
+ c_1 \partial^\mu ( \frac{1}{4} \Sigma^i_\mu{}^\nu \partial_\nu h +2 \partial_\mu h^i -  \Sigma^{i\alpha\beta} \partial_\alpha \tilde{h}_{\mu\beta} ) \\ \nonumber
+  c_2 \epsilon^{ijk} \Sigma^{j\mu\nu} \partial_\mu ( \frac{1}{4} \Sigma^k_\nu{}^\alpha \partial_\alpha h +2 \partial_\nu h^k -  \Sigma^{k\rho\sigma} \partial_\rho \tilde{h}_{\nu\sigma}) = 
-(\sqrt{2} f' +c_1) \Sigma^{i\mu\nu} \partial_\mu \partial^\alpha \tilde{h}_{\nu\alpha}  + 2 c_1 \partial^\mu \partial_\mu h^i. 
\ee
The absence of the first term demands $f'=- c_1/\sqrt{2}$.

We also have, for the operator $(\tilde{d}_1 \tilde{d}_4 + d_2^* \tilde{d}_3^*)\chi$, for the trace part
\be
=  - \sqrt{2} f\Sigma_\mu^i{}^\alpha \partial^\mu \partial_\alpha \chi^i    -2 \Sigma^{i\mu\nu} \partial_\mu (c_1' \partial_\nu \chi^i + c_2' \epsilon^{ijk} \Sigma^j_\nu{}^\alpha \partial_\alpha \chi^k)=0.
\ee
For the vector part
\be
\frac{1}{4} \epsilon^{ijk} \Sigma^{j\mu\nu} \partial_\mu (c_1' \partial_\nu \chi^k + c_2' \epsilon^{klm} \Sigma^l_\nu{}^\alpha \partial_\alpha \chi^m)= \frac{c_2'}{2} \partial^\alpha \partial_\alpha \chi^i.
\ee
For the tracefree part
\be\nonumber
\sqrt{2} \partial_{\langle\mu} f\Sigma_{\nu\rangle}^i{}^\alpha \partial_\alpha \chi^i + \Sigma^i_{\langle\mu}{}^\alpha \partial_{\nu\rangle}(c_1' \partial_\alpha \chi^i + c_2' \epsilon^{ijk} \Sigma^j_\alpha{}^\beta \partial_\beta \chi^k)
-\Sigma^i_{\langle\mu}{}^\alpha \partial_\alpha (c_1' \partial_{\nu\rangle} \chi^i + c_2' \epsilon^{ijk} \Sigma^j_{\nu\rangle}{}^\beta \partial_\beta \chi^k)=\\ \nonumber
(\sqrt{2} f+2c_2') \Sigma^i_{\langle\mu}{}^\alpha \partial_{\nu\rangle} \partial_\alpha \chi^i.
\ee
The vanishing of this requires $f=-\sqrt{2} c_2'$. Note that this is compatible with $2 c_2' = c_1$ and $f=f'$. This fixes all of the operators, apart from the operator $\tilde{d}_3 a$ which can still be chosen arbitrarily.  

We also need to compute
\be
(\tilde{d}_4^* \tilde{d}_4 + \tilde{d}_3 \tilde{d}_3^*) \chi = f' \Sigma^{i\mu\nu} \partial_\mu ( f \Sigma^j_\nu{}^\alpha \partial_\alpha \chi^j) + c_1 \partial^\mu ( c_1' \partial_\mu \chi^i + c_2' \epsilon^{ijk} \Sigma^j_\mu{}^\alpha \partial_\alpha \chi^k) \\ \nonumber
+ c_2 \epsilon^{ijk} \Sigma^{j\mu\nu} \partial_\mu ( c_1' \partial_\nu \chi^k + c_2' \epsilon^{klm} \Sigma^l_\nu{}^\alpha \partial_\alpha \chi^m)
= (-ff' +c_1 c_1' + 2c_2 c_2') \partial^\mu \partial_\mu \chi^i.
\ee

\subsection{The operator \texorpdfstring{$\tilde{D}^* \tilde{D}$}{}}

With the choices we have made, the operator $\tilde{D}^* \tilde{D}$ takes the following form
\be
\tilde{D}^* \tilde{D} \left( \begin{array}{c} h \\ h^i \\ \tilde{h}_{\mu\nu} \\ \chi^i \end{array}\right) = - \partial^\mu \partial_\mu \left( \begin{array}{c} - h \\ - \frac{c_1}{4} \chi^i \\ \tilde{h}_{\mu\nu} \\ -2c_1 h^i + (c_1^2-2c_1 c_2)\chi^i\end{array}\right).
\ee
We have parametrised all the appearing coefficients in terms of $c_1, c_2$. We see that the operator is not diagonal in the space $h^i,\chi^i$, but rather a matrix
\be
\left( \begin{array}{cc} 0 & - \frac{c_1}{4} \\ -2c_1 & c_1^2-2c_1 c_2 \end{array}\right)
\ee
appears. We would like to have the eigenvalues of this matrix to be $\pm 1$, which requires the diagonal elements to be both zero, and the product of the off-diagonal elements to be unity. This means
\be
c_1= \sqrt{2}, \qquad c_2 = \frac{1}{\sqrt{2}}.
\ee
Note that this means $f=-1, c_1'= 0, c_2'= 1/\sqrt{2}$. With these choices the matrix mixing $h^i, \chi^i$ becomes
\be
\left( \begin{array}{cc} 0 & - \frac{1}{2\sqrt{2}} \\ -2\sqrt{2} & 0 \end{array}\right)
\ee

\subsection{Final result for the operator $\tilde{D}$}
We now note that, with the choices we have made, the twisted operator $\tilde{D}$, written as a matrix acting on $S \oplus E$ in terms of the original Pleba\'nski maps becomes 
\be
\tilde{D} = \begin{pmatrix}\frac{1}{\sqrt{2}}(d^*_1 - \Phi d_2) && -\Phi d^*_3 \\ d_2 && \frac{1}{\sqrt{2}} J_1 d^*_3 \end{pmatrix}.
\ee
This should be contrasted with the "naive" map (\ref{D}), which in the same notations is represented by the following matrix
\be
D = \begin{pmatrix}d^*_1 && 0 \\ d_2 &&  d^*_3 \end{pmatrix}.
\ee

\subsection{Lagrangian}

We can obtain the operators $\tilde{D}, \tilde{D}^*$ as arising from an action
\be
S = \int  \xi^\mu  (\tilde{d}_1^* \sigma+ \tilde{d}_4 \chi)_\mu + \epsilon^{ijk} \Sigma^{i\mu\nu} a^j_\mu( d_2\sigma+ \tilde{d}_3^* \chi)_\nu^k - \frac{1}{2} (\xi_\mu)^2 - \frac{1}{2} \epsilon^{ijk} \Sigma^{i\mu\nu} a^j_\mu a^k_\nu.
\ee
Explicitly, substituting the operators, we get
\be
S = \int  \xi^\mu  (\sqrt{2}( \frac{1}{4} \partial_\mu h -  \partial^\nu \tilde{h}_{\mu\nu}) - \Sigma^i_\mu{}^\alpha \partial_\alpha \chi^i) \\ \nonumber
+ \epsilon^{ijk} \Sigma^{i\mu\nu} a^j_\mu( \frac{1}{4} \Sigma^k_\nu{}^\alpha \partial_\alpha h +2 \partial_\nu h^k -  \Sigma^{k\alpha\beta} \partial_\alpha \tilde{h}_{\nu\beta} + \frac{1}{\sqrt{2}} \epsilon^{klm} \Sigma^l_\nu{}^\alpha \partial_\alpha \chi^m) \\ \nonumber
- \frac{1}{2} (\xi_\mu)^2 - \frac{1}{2} \epsilon^{ijk} \Sigma^{i\mu\nu} a^j_\mu a^k_\nu.
\ee
This is very similar to the Lagrangian that appears as (8.174) with (8.178) in \cite{Krasnov:2020lku}, apart from the last term in the second line, which would need to be chosen differently in order to get the Lagrangian in \cite{Krasnov:2020lku}. We have thus reproduced and developed further the gauge-fixing procedure described in this reference. 

Varying this with respect to $\xi, a$ we get
\be
\xi = \tilde{d}_1^* \sigma+ \tilde{d}_4 \chi, \qquad a = d_2\sigma+ \tilde{d}_3^* \chi. 
\ee
Substituting this back into the action we get a second order in derivatives action functional that depends only on $\sigma, \chi$
\be\nonumber
S= \frac{1}{2} \int \left(\sqrt{2}( \frac{1}{4} \partial_\mu h -  \partial^\nu \tilde{h}_{\mu\nu})- \Sigma^i_\mu{}^\alpha \partial_\alpha \chi^i \right)^2 \\ \nonumber
+ \epsilon^{ijk} \Sigma^{i\mu\nu} ( \frac{1}{4} \Sigma^j_\mu{}^\alpha \partial_\alpha h +2 \partial_\mu h^j -  \Sigma^{j\alpha\beta} \partial_\alpha \tilde{h}_{\mu\beta} + \frac{1}{\sqrt{2}} \epsilon^{jlm} \Sigma^l_\mu{}^\alpha \partial_\alpha \chi^m)\times \\ \nonumber
( \frac{1}{4} \Sigma^k_\nu{}^\rho \partial_\rho h +2 \partial_\nu h^k -  \Sigma^{k\rho\sigma} \partial_\rho \tilde{h}_{\nu\sigma} + \frac{1}{\sqrt{2}} \epsilon^{kpq} \Sigma^p_\nu{}^\rho \partial_\rho \chi^q).
\ee
Evaluating this gives
\be\label{sec-order-action}
S=\frac{1}{2} \int - \frac{1}{4} (\partial_\mu h)^2 - 4\sqrt{2} \partial_\mu h^i \partial^\mu \chi^i + ( \partial_\alpha \tilde{h}_{\mu\nu})^2.
\ee
This computation of the second order action is just another way to state that the introduced operator $\tilde{D}$ has the property $\tilde{D}^* \tilde{D}\sim \Delta$.

\section{Splitting the elliptic operator \texorpdfstring{$\tilde{D}$}{}}

We now want to understand the arising elliptic operator $\tilde{D}$ better. 

\subsection{Changing the parametrisation of the \texorpdfstring{$(\xi,a)$}{} space}

Let us introduce new coordinates on the space $E\oplus E\times\Lambda^1$
\be
\Omega_\mu^i = a_\mu^i - \frac{1}{\sqrt{2}} \Sigma^i_\mu{}^\alpha \xi_\alpha, \qquad \omega_\mu = \xi_\mu + \sqrt{2} \Sigma^i_\mu{}^\alpha a_\alpha^i.
\ee
The inverse transformation is
\be
\xi_\mu = \frac{1}{2}\omega_\mu + \frac{1}{\sqrt{2}} \Sigma_\mu^i{}^\alpha \Omega^i_\alpha, \qquad a_\mu^i = \frac{1}{2}( \Omega_\mu^i + \epsilon^{ijk}\Sigma^j_\mu{}^\alpha \Omega_\alpha^k) - \frac{1}{2\sqrt{2}} \Sigma^i_\mu{}^\alpha \omega_\alpha.
\ee
Using the previously introduced operators $\Phi,\Phi^*$ we have
\be
\Omega^i = a^i - \sqrt{2} \Phi^*(\xi)^i, \quad \omega = \xi + \sqrt{2} \Phi(a)
\ee
\be
\xi = \frac{1}{2}\omega + \frac{1}{\sqrt{2}} \Phi(\Omega), \quad a^i = \frac{1}{2}(\Omega^i + J_1(\Omega)^i) - \frac{1}{\sqrt{2}} \Phi^*(\omega).
\ee
The transformation is selected so that
\be\label{inner-omega}
\epsilon^{ijk} \Sigma^{i\mu\nu} a_\mu^j a_\nu^k + (\xi_\mu)^2 = (\Omega_\mu^i)^2 - \frac{1}{2} (\omega_\mu)^2. 
\ee
We have thus found the variables that diagonalise the indefinite quadratic form that appears on the left-hand side of (\ref{inner-omega}). This makes it clear that it is more interesting to consider the operator $\tilde{D}$ composed with the transformation to $\omega, \Omega$ variables. 

\subsection{\texorpdfstring{Splitting of $\tilde{D}$}{}}

Let $T_2$ be the linear map on $(E,E\times\Lambda^1)$ giving $(\omega,\Omega)$ from $(\xi,a)$
\be
\left( \begin{array}{c} \omega \\ \Omega \end{array}\right) = \left( \begin{array}{cc} \mathbb{I} & \sqrt{2}\Phi \\ - \sqrt{2} \Phi^* & \mathbb{I}\end{array}\right) \left( \begin{array}{c} \xi \\ a \end{array}\right).
\ee
We then have
\be
T_2 \tilde{D} = \begin{pmatrix} 1 && \sqrt{2} \Phi \\ - \sqrt{2} \Phi^* && 1 \end{pmatrix} \begin{pmatrix} \frac{1}{\sqrt{2}}(d_1^* - \Phi d_2) && - \Phi d_3^* \\ d_2 && \frac{1}{\sqrt{2}} J_1 d_3^* \end{pmatrix} =\\ \nonumber \begin{pmatrix} \frac{1}{\sqrt{2}} d_1^* + (\sqrt{2}-\frac{1}{\sqrt{2}})  \Phi d_2  && - \Phi d_3^* + \Phi J_1 d_3^* \\ d_2 - \Phi^*(d_1^* - \Phi d_2) && \sqrt{2} \Phi^* \Phi d_3^* + \frac{1}{\sqrt{2}} J_1 d_3^* \end{pmatrix}
\ee
To simplify this further we will need the following identities 
\be
(\Phi J_1)(a)_\mu = \Sigma^i_\mu{}^\nu \epsilon^{ijk} \Sigma^j_\nu{}^\rho a^k_\rho = 2 \Sigma^i_\mu{}^\nu a^i_\nu = 2 \Phi(a)_\mu,
\ee
\be
(\Phi^* \Phi)(a)^i_\mu = \frac{1}{2} \Sigma^i_\mu{}^\nu \Sigma^j_\nu{}^\rho a^j_\rho = -\frac{1}{2} a^i_\mu - \frac{1}{2} \epsilon^{ijk}\Sigma^j_\mu{}^\nu a^k_\nu = -\frac{1}{2}(1+J_1)(a)^i_\mu.
\ee
Using the language of operators the above is simply 
\be
\Phi J_1 = 2 \Phi, \qquad \Phi^* \Phi = -\frac{1}{2}(1+J_1).
\ee
This gives
\be
T_2 \tilde{D} =\begin{pmatrix} \frac{1}{\sqrt{2}} (d_1^* +  \Phi d_2)  &&  \Phi d_3^*  \\ \frac{1}{2}(1-J_1) d_2 - \Phi^*d_1^*  &&- \frac{1}{\sqrt{2}}  d_3^* \end{pmatrix}
\ee
The operators appearing in the first column can be further simplified. We have
\be
(d_1^* +  \Phi d_2)\sigma =  -\frac{1}{4}\partial_\mu h + 2\Sigma^i_\mu{}^\nu \partial_\nu h^i -  \partial^\nu \tilde{h}_{\mu\nu}  - \frac{3}{4} \partial_\mu h + \partial^\nu \tilde{h}_{\mu\nu} + 2 \Sigma^i_\mu{}^\alpha\partial_\alpha h^i = \\ \nonumber
- \partial_\mu h + 4\Sigma^i_\mu{}^\nu \partial_\nu h^i,
\ee
where we took $\Phi d_2$ from (\ref{Phi-d2}). Note cancellation of the $\tilde{h}_{\mu\nu}$ terms. We also have
\be
\Phi^* d_1^* = \frac{1}{2} \Sigma^i_\mu{}^\alpha ( -\frac{1}{4}\partial_\alpha h + 2\Sigma^j_\alpha{}^\nu \partial_\nu h^j -  \partial^\nu \tilde{h}_{\alpha\nu} ) = \\ \nonumber
- \frac{1}{8} \Sigma^i_\mu{}^\alpha \partial_\alpha h - \partial_\mu h^i - \epsilon^{ijk} \Sigma^j_\mu{}^\alpha \partial_\alpha h^k - \frac{1}{2}  \Sigma^i_\mu{}^\alpha \partial^\nu \tilde{h}_{\alpha\nu} , 
\ee
and
\be
J_1 d_2 \sigma = \epsilon^{ijk} \Sigma^j_\mu{}^\alpha ( \frac{1}{4} \Sigma^k_\alpha{}^\nu \partial_\nu h +2 \partial_\alpha h^k -  \Sigma^{k\rho\sigma} \partial_\rho \tilde{h}_{\alpha\sigma} )= \\ \nonumber
\frac{1}{2} \Sigma^i_\mu{}^\nu \partial_\nu h + 2 \epsilon^{ijk} \Sigma^j_\mu{}^\alpha \partial_\alpha h^k + \Sigma^{i\alpha\beta} \partial_\alpha \tilde{h}_{\beta\mu}-
\Sigma^{i}_\mu{}^\alpha \partial^\beta \tilde{h}_{\alpha\beta},
\ee
so that
\be
(1-J_1)d_2\sigma =- \frac{1}{4} \Sigma^i_\mu{}^\nu \partial_\nu h +2 \partial_\mu h^i  - 2 \epsilon^{ijk} \Sigma^j_\mu{}^\alpha \partial_\alpha h^k - 2\Sigma^{i\alpha\beta} \partial_\alpha \tilde{h}_{\beta\mu}+
\Sigma^{i}_\mu{}^\alpha \partial^\beta \tilde{h}_{\alpha\beta}.
\ee
This gives, finally
\be\nonumber
\left(\frac{1}{2}(1-J_1) d_2 - \Phi^*d_1^*\right)\sigma = - \frac{1}{8} \Sigma^i_\mu{}^\nu \partial_\nu h + \partial_\mu h^i  -  \epsilon^{ijk} \Sigma^j_\mu{}^\alpha \partial_\alpha h^k - \Sigma^{i\alpha\beta} \partial_\alpha \tilde{h}_{\beta\mu}+\frac{1}{2} 
\Sigma^{i}_\mu{}^\alpha \partial^\beta \tilde{h}_{\alpha\beta} \\ \nonumber
+ \frac{1}{8} \Sigma^i_\mu{}^\alpha \partial_\alpha h + \partial_\mu h^i + \epsilon^{ijk} \Sigma^j_\mu{}^\alpha \partial_\alpha h^k + \frac{1}{2}  \Sigma^i_\mu{}^\alpha \partial^\nu \tilde{h}_{\alpha\nu} = \\ \nonumber
2 \partial_\mu h^i+\Sigma^{i}_\mu{}^\alpha \partial^\beta \tilde{h}_{\alpha\beta}- \Sigma^{i\alpha\beta} \partial_\alpha \tilde{h}_{\beta\mu}.
\ee
Note cancelation of the $h$-dependent terms here. 

We now rewrite the operator $T_2 \tilde{D}$ as one acting on a multiplet of fields $(h,h^i,\chi^i,\tilde{h}_{\mu\nu})$. Using matrix notation we have
\be
\left( \begin{array}{c} \omega \\ \Omega\end{array}\right)=
T_2 \tilde{D} \left( \begin{array}{c} h \\ h^i \\ \chi^i \\ \tilde{h}_{\mu\nu}\end{array}\right) = \left(\begin{array}{cccc} - \frac{1}{\sqrt{2}} d & 2\sqrt{2} \Phi d & \Phi d & 0 \\
0 & 2 d & - \frac{1}{\sqrt{2}} d & \tilde{d} \end{array}\right) \left( \begin{array}{c} h \\ h^i \\ \chi^i \\ \tilde{h}_{\mu\nu}\end{array}\right),
\ee
where we have defined new operators 
\be
d:\Lambda^0\to \Lambda^1, \qquad & (dh)_\mu := \partial_\mu h, \\ \nonumber
d: E\to E\otimes \Lambda^1, \qquad & (dh^i)_\mu := \partial_\mu h^i, \\ \nonumber
\tilde{d}:{\rm Sym}^2_0(\Lambda^1)\to E\otimes \Lambda^1, \qquad & (\tilde{d}\tilde{h})_\mu^i := \Sigma^{i}_\mu{}^\alpha \partial^\beta \tilde{h}_{\alpha\beta}- \Sigma^{i\alpha\beta} \partial_\alpha \tilde{h}_{\beta\mu}.
\ee
We note that we can rewrite the operator $\tilde{d}$ in a different way using the $\diamond$ operator
\be
(\tilde{h}\diamond \Sigma^i)_{\mu\nu}:= \tilde{h}_\mu{}^\alpha \Sigma^i_{\alpha\nu} -  \tilde{h}_\nu{}^\alpha \Sigma^i_{\alpha\mu}.
\ee
Then
\be
 (\tilde{d}\tilde{h})_\nu^i = - \partial^\mu (\tilde{h}\diamond \Sigma^i)_{\mu\nu}.
 \ee
Let also write the arising objects $\omega,\Omega$ explicitly
\be
\omega_\mu = - \frac{1}{\sqrt{2}} \partial_\mu h  + \Sigma^i_\mu{}^\alpha \partial_\alpha ( 2\sqrt{2} h^i + \chi^i), \\ \nonumber
\Omega_\mu^i =\Sigma^i_\mu{}^\alpha\partial^\beta \tilde{h}_{\alpha\beta} -  \Sigma^{i\alpha\beta} \partial_\alpha \tilde{h}_{\mu\beta} + \frac{1}{\sqrt{2}} \partial_\mu (2 \sqrt{2} h^i - \chi^i).
\ee

We now introduce
\be
h_\pm^i := 2 h^i \pm \frac{1}{\sqrt{2}} \chi^i.
\ee
Note that these are precisely the variables that diagonalise the mixed term in (\ref{sec-order-action}). 
\be
4 \sqrt{2} \partial^\mu h^i \partial_\mu \chi_i =  ( ( \partial_\mu h_+)^2 -  ( \partial_\mu h_-)^2).
\ee
We also have
\be
8 (h^i)^2 + (\chi^i)^2 = (h_+^i)^2 + (h_-^i)^2,
\ee
so that these are also the correct variables to diagonalise the inner product quadratic form. We can then write
\be
\left( \begin{array}{c} h\\ h^i \\ \chi^i \\ \tilde{h}_{\mu\nu} \end{array}\right) = 
 \left( \begin{array}{cccc}  1& 0 & 0 & 0 \\ 0 & \frac{1}{4} & \frac{1}{4} & 0 \\ 0&  \frac{1}{\sqrt{2}} & - \frac{1}{\sqrt{2}} &0 \\ 0 & 0 & 0 & 1 \end{array} \right) \left( \begin{array}{c} h\\ h^i_+ \\ h^i_- \\ \tilde{h}_{\mu\nu} \end{array} \right) = T_1  \left( \begin{array}{c} h\\ h^i_+ \\ h^i_- \\ \tilde{h}_{\mu\nu} \end{array} \right) .
\ee
We then have
\be
T_2 \tilde{D} T_1 = \left(\begin{array}{cccc} - \frac{1}{\sqrt{2}} d & \sqrt{2} \Phi d & 0 & 0 \\
0 & 0 &d & \tilde{d} \end{array}\right).
\ee
This shows that the operator $T_2 \tilde{D} T_1$ splits into the direct sum of two operators. 

Let us describe the arising operators more explicitly. We define
\be
D_4: \Lambda^0 \oplus E \to \Lambda^1, \qquad D_{12}: E \oplus {\rm Sym}^2_0(\Lambda^1) \to E\otimes \Lambda^1,
\ee
given by
\be
D_4 (h,h_+) = \left( \begin{array}{cc} - \frac{1}{\sqrt{2}} d & \sqrt{2} \Phi d\end{array}\right) \left( \begin{array}{c} h \\ h_+\end{array}\right), \\ \nonumber
D_{12} (h_-, \tilde{h}) = \left( \begin{array}{cc} d & \tilde{d} \end{array}\right) \left( \begin{array}{c} h_- \\ \tilde{h} \end{array}\right).
\ee
Explicitly, these operators are given by
\be\label{d4-d12}
D_4 (h,h^i_+) = 
- \frac{1}{\sqrt{2}} \left( \partial_\mu h - 2\Sigma^i_\mu{}^\alpha \partial_\alpha h^i_+ \right), \\ \nonumber
 D_{12} (h^i_-, \tilde{h}_{\mu\nu})  =  \partial_\mu h^i_- + \Sigma^i_\mu{}^\alpha\partial^\beta \tilde{h}_{\alpha\beta} -  \Sigma^{i\alpha\beta} \partial_\alpha \tilde{h}_{\mu\beta}.
\ee

Each of these two operators times its adjoint gives a multiple of $\Delta$. This is easiest confirmed as an action computation. We have 
\be
\int (D_4 (h,h_+) )^2 = \frac{1}{2} \int  (\partial_\mu h)^2 + 4(\partial_\alpha h^i_+)^2, \\ \nonumber
\int (D_{12} (h_-, \tilde{h}))^2 = \int (\partial_\alpha h^i_-)^2 + (\partial_\alpha \tilde{h}_{\mu\nu})^2.
\ee
Here we used
\be\nonumber
\int ( \Sigma^i_\mu{}^\alpha\partial^\beta \tilde{h}_{\alpha\beta} -  \Sigma^{i\alpha\beta} \partial_\alpha \tilde{h}_{\mu\beta})^2=\int 3 (\partial^\nu \tilde{h}_{\mu\nu})^2 - 2 (\partial^\nu \tilde{h}_{\mu\nu})^2+ (\partial_\alpha \tilde{h}_{\mu\nu})^2 - (\partial^\nu \tilde{h}_{\mu\nu})^2 = \int (\partial_\alpha \tilde{h}_{\mu\nu})^2.
\ee
We then have
\be
\int (D_{12} (h_-, \tilde{h}))^2 - \frac{1}{2} (D_4 (h,h_+) )^2 = \int (\partial_\alpha \tilde{h}_{\mu\nu})^2 + (\partial_\alpha h^i_-)^2 - (\partial_\alpha h^i_+)^2 - \frac{1}{4} (\partial_\mu h)^2 = 2S,
\ee
where $S$ is the action (\ref{sec-order-action}). 


\end{document}